\documentclass[dvips,ps]{imsart}

\usepackage{amsthm}
\usepackage{amsmath}
\usepackage[psamsfonts]{amssymb}
\RequirePackage[numbers]{natbib}
\RequirePackage[colorlinks,citecolor=blue,urlcolor=blue]{hyperref}
\usepackage{graphics}
\doi{10.1214/154957804100000000}
\pubyear{0000}
\volume{0}
\firstpage{0}
\lastpage{0}

\startlocaldefs
\newcommand{\mc}{\mathcal}
\newcommand{\mf}{\mathfrak}
\newcommand{\dg}{\dagger}
\newcommand{\ind}{{\bf 1}}
\newcommand{\Z}{\mathbb{Z}}
\newcommand{\R}{\mathbb{R}}
\DeclareMathOperator{\Hom}{Hom}
\newcommand{\C}{\mathbb{C}}
\newcommand{\U}{\mathbb{U}}
\newcommand{\E}{\mathbb{E}}
\newcommand{\N}{\mathbb{N}}
\DeclareMathOperator{\Id}{Id}
\DeclareMathOperator{\End}{End}
\DeclareMathOperator{\Ker}{Ker}
\DeclareMathOperator{\Cov}{Cov}
\newcommand{\eps}{\varepsilon}
\newcommand{\Lap}{\Delta\!}
\newcommand{\dmd}{\diamondsuit}
\renewcommand{\H}{\mathbb{H}}
\renewcommand{\P}{\mathbb{P}}

\endlocaldefs

\begin{document}

\begin{frontmatter}

\title{Topics on abelian spin models and related problems}
\thankstext{t1}{Partially supported by NSF grant DMS-1005749 and the Alfred P. Sloan Foundation.}

\runtitle{Topics on abelian spin models and related problems}

\author{\fnms{Julien} \snm{Dub\'edat\thanksref{t1}}\ead[label=e1]{dubedat@math.columbia.edu}}
\address{
Columbia University, Department of Mathematics\\
2990 Broadway, New York, NY 10027\\
 \printead{e1}}

\runauthor{J. Dub\'edat}

\begin{abstract}
In these notes, we discuss a selection of topics on several models of planar statistical mechanics. We consider the Ising, Potts, and more generally abelian spin models; the discrete Gaussian free field; the random cluster model; and the six-vertex model. Emphasis is put on duality, order, disorder and spinor variables, and on mappings between these models.
\end{abstract}


\end{frontmatter}

\section{Introduction}

In the type of statistical models considered here, a fixed underlying graph carries a random realization of a structure, consisting of values (boolean, in a finite alphabet or scalar) carried by bonds (edges) or sites (vertices), and usually interacting via local rules. We will focus on the Ising and Potts models (and more generally cyclic and abelian spin models), where each site carries a random ``spin"; the discrete Gaussian free field, where a scalar height function is defined on sites; the random cluster model, where realizations are random subgraphs; and the six-vertex model, where each edge has a random orientation.

The key question is to understand the large scale random structures  generated by these local interactions. For instance, one may study the correlations between local configurations observed at large distance on the underlying graph. We will generically refer to (functions of) these local inputs as order variables.

Historically, a fundamental tool, introduced by Kramers and Wannier in the context of the Ising model, is duality, which maps the model defined on a (weighted) graph to the same model on another (weighted) graph. Of particular interest are the ``fixed points" of these duality transformations; they often coincide with a critical point of the model, at which one can observe large scale random structures. This type of duality is observed in a variety of models, and for abelian models can be phrased in terms of the Fourier-Pontryagin transform and the Poisson summation formula. 

Again in the context of the Ising model, Kadanoff and Ceva showed that under duality, order variables are mapped to ``disorder" variables, which represent a modification of the state space by the introduction of local ``defects"; moreover, the field-theoretic concept of fermion finds a combinatorial incarnation as the combination of an order and a disorder variable at microscopic distance. Such combinations are referred to as parafermionic or spinor variables.

In these notes, we will describe duality, order, disorder, and spinor variables for abelian spin models and the discrete free field. We then discuss combinatorial mappings between the Potts model, the random cluster model, the six-vertex model, and the dimer model.

For focus and (relative) self-containedness, we do not discuss here a number of major contiguous subjects. These include: phase transition, transfer matrices and Bethe ansatz (\cite{Baxterexact,Resh6V}), Yang-Baxter solvability (\cite{Baxterexact}), infinite volume measures (\cite{GeoGibbs,GrimFK}) and limit shapes (\cite{CKP,KenIAS,Resh6V}), scaling limits and Schramm-Loewner Evolution (\cite{W1,SmiICM}).

\section{Spin models}

\subsection{Ising model}

{\bf Definition}\\

In the Ising model (\cite{MWising,Palmerplanar,Grimgraphs}) , a configuration on an underlying finite graph $\Gamma=(V,E)$ consists of an assignment of a spin value $\sigma(v)\in\{1,-1\}$ to each vertex $v\in V$. The weight of a configuration is:
$$w(\sigma)=\prod_{e=(vv')\in E}\exp\left(\beta J_e\sigma(v)\sigma(v')\right)$$
where $J_e\neq 0$ is a coupling constant attached to the edge $e\in E$ and $\beta>0$ is the inverse temperature. 
Plainly the model (up to normalization) depends only on the edge weights:
$$w(e)=\exp(-2\beta J_e)$$
It is occasionally convenient to allow $w(e)=0$, i.e. $J_e=+\infty$. Note that these weights are invariant under spin flip: $w(-\sigma)=w(\sigma)$. One may also consider $\pm$ boundary conditions, that is, fixing the values of spins to $\pm 1$ on some prescribed subsets $V^\pm$ of $V$ (the configuration space then consists of spin collections $(\sigma(v))_{v\in V}$ s.t. $\sigma(v)=\pm 1$ for $v\in V^\pm$). One may also consider some segments on the boundary to be {\em wired}, viz. connected by edges with zero weight, so that the only contributing configurations are constant on these segments (a wired arc may be equivalently represented as an extended vertex). 

The {\em partition function} is 
$${\mc Z}={\mc Z}(\Gamma,(J_e)_{e\in E})=\sum_{\sigma:V\rightarrow\{\pm1\}}w(\sigma)$$
We denote the expectation under the probability measure given by $\P(\sigma=\sigma_0)=w(\sigma_0)/{\mc Z}$ by $\langle\cdot\rangle$. 

A classical object of interest is the {\em spin correlation}
$$\langle \sigma(v_1)\cdots\sigma(v_n)\rangle$$
for $v_1,\dots,v_n$ graph vertices (in particular when they are at large distance) and the spin variables $\sigma(v)$'s are fundamental examples of {\em order variables}. 

The spin correlations depend on the boundary conditions. A wired arc may be seen as a boundary arc along which couplings are infinite (vanishing edge weights); in turn, $\pm$ boundary arcs may be represented in terms of wired arcs and spin variables. For example, if $\gamma^+,\gamma^-$ are two disjoint boundary arcs, we have
$$\langle X\rangle_{+,-}=\frac{\left\langle X\cdot\frac{1+\sigma(x)}2\cdot\frac{1-\sigma(y)}2\right\rangle_{w,w}}{\left\langle\frac{1+\sigma(x)}2\cdot\frac{1-\sigma(y)}2\right\rangle_{w,w}}$$
where $X$ is a generic random variable, the expectation on the LHS is for $\pm$ boundary conditions on $\gamma^\pm$, and the expectations on the RHS are for $\gamma^+$, $\gamma^-$ (separately) wired arcs. 
This shows that, allowing for general couplings, it is enough to treat free boundary conditions.

\begin{figure}[htb]
\begin{center}
\scalebox{.4}{\includegraphics{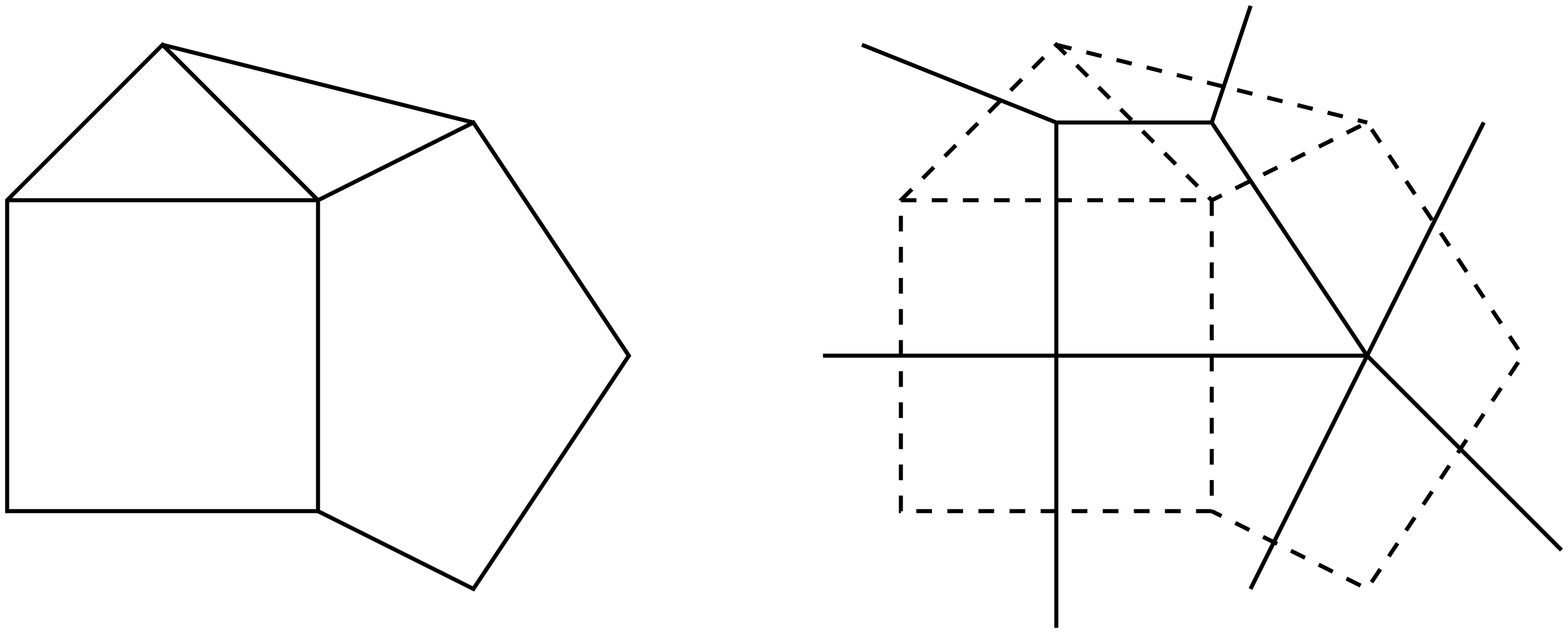}}
\end{center}
\caption{Left: a portion of a planar graph $\Gamma$. Right: its dual $\Gamma^\dg$}
\label{Fig:dual}
\end{figure}

From now on we will consider the planar case; $\Gamma^\dg=(V^\dg,E^\dg)$ denotes the planar dual of $\Gamma$, so that faces of $\Gamma$ correspond to vertices of $\Gamma^\dg$, and vice-versa (Figure \ref{Fig:dual}). 

More precisely, it is rather convenient to think of vertices on a wired arc of the boundary as a single extended vertex; and to connect vertices on a free boundary or extended vertices  to a vertex at infinity by edges with weight 1 (i.e. with no interaction). If $e\in E$, we denote $e^\dg$ the corresponding edge of $\Gamma^\dg$; for oriented edges, $e^\dg$ crosses $e$ from right to left.\\

{\bf Graphical expansions}\\

In the low-temperature regime ($\beta\gg 1$), disagreements between neighboring spins are severely penalized. The {\em low temperature expansion} of the Ising model consists in mapping a spin configuration $(\sigma_v)_{v\in V}$ to a subgraph $P=(V^\dg,E_P)$ of $\Gamma^\dg$, where $E_P=\{e=(vv')^\dg\in E^\dg: \sigma(v)\sigma(v')=-1\}$. At low temperature this graph is typically sparse. The spins $(\sigma_v)$ may be recovered from $P$ up to an overall spin flip (in the case of free boundary conditions). It is easily checked that $P$ is an even degree subgraph of $P$; such subgraphs are called {\em polygons}. 

The (unnormalized) Ising measure on spins $(\sigma_v)$ projects to a measure on polygons:
$$w(P)=\ind_{\{P{\rm\ admissible}\}}\prod_{e\in E_P}w(e^\dg)$$
omitting a multiplicative factor $2$ (for overall spin flip).

Let us turn to the {\em high temperature expansion}. Since $\sigma(v)\sigma(v')=\pm 1$, we may write
$$\exp(\beta J_e\sigma(v)\sigma(v'))=\cosh(\beta J_e)\left (1+\tanh(\beta J_e)\sigma(v)\sigma(v')\right)$$
Let us define a dual edge weight:
$$w'(e)=\tanh(\beta J_e)$$
In the high temperature regime ($\beta\ll 1$), these dual edge weights are close to 1. Starting from the partition function, we obtain
\begin{align*}
{\mc Z}(\Gamma,(J_e)_e)&=\sum_{\sigma:V\rightarrow\{\pm1\}}\prod_{e=(vv')\in E}\exp(-2\beta J_e\ind_{\{\sigma(v)\sigma(v')=-1\}})\\
&=\prod_{e\in E}\cosh(\beta J_e)\sum_{\sigma:V\rightarrow\{\pm1\}}\prod_{e=(vv')\in E}\left (1+w'(e)\sigma(v)\sigma(v')\right)\\
&=\prod_{e\in E}\cosh(\beta J_e)\sum_{\sigma:V\rightarrow\{\pm1\}}\sum_{E_0\subset E}\prod_{e\in E_0}w'(e)\sigma(v)\sigma(v')
\end{align*}
by fulling expanding the product. Then exchanging summations leads to
$${\mc Z}(\Gamma,(w_e)_e)=2^{|V|}\prod_{e\in E}\cosh(\beta J_e)\sum_{P{\rm\ polygon}}\prod_{e\in E_P}w'(e)$$
Then comparing high and low temperature expansions, we have the {\em Kramers-Wannier duality} (\cite{KWising1,KWising2}) for partition functions:
$${\mc Z}(\Gamma,(J_e)_e)=2^{|V|}\left(\prod_{e\in E}\cosh(\beta J_e)\right){\mc Z}(\Gamma^\dg,(J^\dg_e))$$
where
$$e^{-2\beta J^\dg_e}=\tanh(\beta J_e)$$
or $w'(e)=\frac{1-w(e)}{1+w(e)}$. Note that $w\mapsto\frac{1-w}{1+w}$ is involutive. It exchanges $0$ and $1$, so that wired boundary components are exchanged with free boundary components. To retain positive weights after the duality transformation, we need to start from ferromagnetic couplings, i.e. $J_e\geq 0$.

If $\Gamma$ is (a portion of) the square lattice, $\Gamma^\dg$ is identical to $\Gamma$, up to boundary modifications. For instance, one can consider $\Gamma$ to be a rectangle with half-wired and half-free boundary conditions, in which case $\Gamma^\dg$ is exactly isomorphic (as a weighted graph) to $\Gamma$. Thus for the square lattice the fixed point $w_{sd}=\sqrt 2-1$ of the duality mapping $w\mapsto\frac{1-w}{1+w}$ is {\em self-dual}, in the sense that it is fixed for Kramers-Wannier duality (at least ``in the bulk"). The low-temperature regime $w<w_{sd}$ and the high-temperature regime $w>w_{sd}$ are exchanged by duality.
\\

{\bf Disorder variables.}\\

Let us now consider the effect of Kramers-Wannier on spin correlations $\langle \sigma(v_1)\cdots\sigma (v_n)\rangle$. In the high temperature expansion (and periodic boundary conditions for simplicity), we can write:
\begin{align*}
{\mc Z}(\Gamma,(J_e)_e)\langle \sigma(v_1)\cdots\sigma (v_n)\rangle&=\sum_{\sigma\in\{\pm 1\}^{V\setminus V^\pm}}\sigma(v_1)\cdots\sigma(v_n)\prod_{e=(vv')\in E}e^{\beta J_e\sigma(v)\sigma(v')}
\end{align*}
and then 
\begin{align*}
\frac{{\mc Z}(\Gamma,(J_e)_e)}{\prod_{e\in E}\cosh(\beta J_e)}
\langle \sigma_{v_1}\cdots\sigma_{v_n}\rangle
&=
\sum_{\sigma\in\{\pm 1\}^{V\setminus V^\pm}}\sum_{E_0\subset E}\sigma_{v_1}\cdots\sigma_{v_n}\prod_{e\in E_0}w'(e)\sigma_v\sigma_{v'}\\
&=2^{|V|}\sum_{P\subset\Gamma}\prod_{e\in E_P}w'(e)
\end{align*}
where the sum bears on subgraphs $P\subset\Gamma$ which have odd degree at $\{v_1,\dots,v_n\}$ and even degree elsewhere (by $\sigma\leftrightarrow -\sigma$ invariance, the correlator is zero if $n$ is odd, in which case the sum is empty). For periodic boundary conditions (i.e. when $\Gamma$ is embedded on a torus), one requires additionally that any non-contractible (with respect to the torus) path on $\Gamma^\dg$ crosses an even number of edges of $P$.

Such a polygon with defects at $\{v_1,\dots,v_n\}$ does not correspond to a spin configuration on $\Gamma^\dg$. Following Kadanoff and Ceva (\cite{KCdisorder}), this motivates the introduction of {\em disorder variables}, which encode a modification of the state space.

The data of $(\sigma_v)_{v\in V}$ is equivalent, up to global spin flip, to the data $(d\sigma(e))_{e\in E}$, where $d\sigma(vv')=\sigma(v')\sigma(v)^{-1}$. If $\omega:E\rightarrow\{\pm 1\}$, we define:
$$d\omega(f)=\prod_{e\in\partial f}\omega(e)$$
where $f$ is a face of $\Gamma$ and $\partial f$ is its (counterclockwise oriented) boundary. Up to spin flip, spin configurations $(\sigma_v)_{v\in V}$ correspond to closed currents: $\{(\omega_e)_{e\in E}: d\omega\equiv 1\}$. (In the periodic case, there is an additional condition: $\prod_{e\in\gamma}\omega(e)=1$ if $\gamma$ is a non-contractible path). Plainly the Gibbs weights can be written as functions of these currents.

Introducing disorder variables $\mu(f_1),\dots,\dots\mu(f_n)$ ($f_1,\dots,f_n$ faces of $\Gamma$) consists in modifying the state space to:
$$\{(\omega_e)_{e\in E}: d\omega=(-1)^{\ind_{\{f_1,\dots,f_n\}}}\}$$
with weights
$$w(\omega)=\prod_{e\in E:\omega_e=-1} w(e)$$
We assume that $n$ is even (otherwise the state space is empty, as $\prod_{f\in\Gamma^\dg}d\omega(f)=1$ for all $\omega$). If $v,v'\in V$, one can define $\sigma(v')\sigma(v)=\prod_{e\in\gamma_{v\rightarrow v'}}\omega(e)$ where $\gamma_{v\rightarrow v'}$ is a path from $v$ to $v'$ on $\Gamma$. Notice however that this definition depends on the choice of $\gamma$. Indeed, if $\gamma'$ is another such path, $\prod_{e\in\gamma}\omega(e)=\pm\prod_{e\in\gamma}\omega(e)$ depending on the parity of the number of $f_i$'s enclosed by the loop $\gamma'\gamma^{-1}$.

One can trivialize this data by choosing $\frac n2$ ``defect lines" (paths on $\Gamma^\dg$) joining the $f_i$'s pairwise. Then one can find $\sigma:V\rightarrow\{\pm 1\}$ such that $\omega(e)=d\sigma(e)$ if $e$ does not cross one the defect lines. In terms of this trivialization, we have an Ising model with antiferromagnetic (negative) couplings $-J_e$ for $e\in E$ crossing a defect line.

We can now define mixed order-disorder correlators:
$$\langle\prod_i\sigma(v_i)\prod_j\mu(f_j)\rangle{\mc Z}(\Gamma)=
2\sum_{\stackrel{\omega:E\rightarrow\{\pm 1\}}{d\omega=(-1)^{\ind_F}}}
\prod_{e\in\gamma}\omega(e)\prod_{e\in E:\omega(e)=-1}w(e)
$$
where $\{v_1,\dots,v_{2m}\}\in\Gamma$, $F=\{f_1,\dots,f_{2n}\}\subset\Gamma^\dg$, $\gamma$ is a union of $m$ paths on $\Gamma$ with endpoints $v_1,\dots,v_{2m}$. The sign of the correlator depends on this implicit choice of paths (or alternatively of defect lines, in which case one simply requires that $\gamma$ does not intersect defect lines).

In this case the Kramers-Wannier duality (\cite{KWising1,KWising2,KCdisorder}) yields:
$$\langle\prod_i\sigma(v_i)\prod_j\mu(f_j)\rangle_{\Gamma}=\langle\prod_i\mu(v_i)\prod_j\sigma(f_j)\rangle_{\Gamma^\dg}$$
where the couplings on $\Gamma^\dg$ are in duality with those on $\Gamma$.

Kadanoff and Ceva (\cite{KCdisorder}) identified the combination of an order variable $\sigma$ and a disorder variable $\mu$ at microscopic distance as a discrete version of the field-theoretic notion of {\em fermion}. In particular one denotes $\psi_{vf}=\sigma_v\mu_f$ where $v\in V$ is on the boundary of the face $f$.

\subsection{Abelian spin models}

The Ising model may be generalized in several directions; the first we consider here is that of abelian spin models, in which the spin variable at each vertex takes values in a fixed finite abelian group $G$. In the case $G=\{\pm 1\}$, one recovers the Ising model.

Given a graph $\Gamma=(V,E)$ and a finite abelian group $G$, a configuration is a mapping $\sigma:V\rightarrow G$. For each $e\in E$, we have a weight function $w_e:G\rightarrow [0,\infty)$, which is symmetric, i.e. $w_e(g^{-1})=w_e(g)$ (relaxing this condition leads to {\em chiral} models and requires the underlying graph $\Gamma$ to be oriented). Then
$$w(\sigma)=\prod_{e=(vv')\in E}w(\sigma(v')\sigma(v)^{-1})$$
is the weight of a general configuration. For free boundary conditions, we have $w(\rho\sigma)=w(\sigma)$ for any $\rho\in G$.

If $G=\Z/q\Z$, we have a {\em clock model}. If furthermore $w_e=a+b\ind_{\{1\}}$, we get the {\em $q$-state Potts model}. If $q=4$ and $w_e$ is a general symmetric weight, we obtain the {\em Ashkin-Teller model} (\cite{Baxterexact}, as phrased by Fan).

Let $\hat G=\Hom (G,\C^*)$ be the dual of $G$ ($\hat G$ and $G$ are isomorphic, though not canonically). The basic order variables are of type $\chi(\sigma(v))$, $\chi\in\hat G$, with associated correlators:
$$\langle\prod_i\chi_i(\sigma(v_i))\rangle$$
where $v_i\in V$, $\chi_i\in \hat G$. For free boundary conditions, this is zero unless $\prod_i\chi_i=1$.

Kramers-Wannier duality has been extended successively to various models, see in particular \cite{WWduality,Savitdual}. By Fourier-Pontryagin duality (\cite{Rudinfourier}), one can write:
$$w_e=\frac{1}{|G|^{1/2}}\sum_{\chi\in\hat G}\widehat{w_e}(\chi)\chi$$
where $\widehat w_e$ is the Fourier-Pontryagin transform
$$\widehat w_e(\chi)=\frac{1}{|G|^{1/2}}\sum_{g\in G}w_e(g)\bar\chi(g)$$
At the level of partition functions, we have the summation:
\begin{align*}
{\mc Z}(\Gamma,G,(w_e)_e)&=\sum_{\sigma:V\rightarrow G}\prod_{e=(vv')\in E}w_e(\sigma(v')\sigma(v)^{-1})\\
&=|G|^{-\frac 12|E|}\sum_{\sigma:V\rightarrow G}\prod_{e=(vv')\in E}\sum_{\chi_e\in \hat G}\widehat{w_e}(\chi_e)\chi_e(\sigma(v')\sigma(v)^{-1})\\
&=|G|^{-\frac 12|E|}\sum_{\stackrel{\sigma:V\rightarrow G}{\chi:E\rightarrow \hat G}}\prod_{e=(vv')}\widehat{w_e}(\chi_e)\chi_e(\sigma(v')\sigma(v)^{-1})\\
&=|G|^{|V|-\frac 12|E|}\sum_{\chi:E\rightarrow \hat G:d\ast\chi\equiv 1}\prod_{e=(vv')}\sum_{\sigma_v\in G}\widehat{w_e}(\chi_e)
\end{align*}
where $d\ast\chi(v)=\prod_{v'\sim v}\chi_e$. Let us identify $\chi:E\rightarrow\hat G$ with $\chi:E^\dg\rightarrow\hat G$ via $\chi(e^\dg)=\chi(e)$ for $e\in E$. Then $d\ast\chi\equiv 1$ is equivalent (in the simply connected case) to the existence of $\hat \sigma:V^\dg\rightarrow\hat G$ such that $\chi(ff')=\hat g(f')\hat g(f)^{-1}$. (Implicitly, we fixed a reference orientation for edges of $\Gamma^\dg$; by symmetry of the weight, one can write each term $w_e(\sigma(v')\sigma(v)^{-1})$ so that $v$ (resp. $v'$) is the left (resp. right) vertex wrt the oriented edge of $\Gamma^\dg$ separating them). We conclude that:
$${\mc Z}(\Gamma,G,(w_e)_e)=|G|^{|V|-\frac 12|E|+1}{\mc Z}(\Gamma^\dg,\hat G,(\widehat{w_e})_e)$$
Disorder variables are indexed by a face $f\in V^\dg$ and an element $g\in G$, and denoted $\mu_g(f)$. The corresponding modified state space consists of $G$-valued 1-forms $\omega:\overrightarrow E\rightarrow G$ defined on oriented edges and antisymmetric in the sense that $\omega(\overrightarrow{xy})=\omega(\overrightarrow{yx})^{-1}$. We define $d\omega(f)=\prod_{\overrightarrow e\in\partial f}\omega(\overrightarrow e)$, where $\partial f$ is the counterclockwise oriented boundary of $f$. Then the state space corresponding to a disorder $\prod_i\mu_{g_i}(f_i)$ is:
$$\{\omega: d\omega=\prod_i g_i^{\ind_{f_i}}\}$$
with weight $w(\omega)=\prod_{e\in E}w_e(\omega(e))$ (there is an overall factor $|G|$ compared with the earlier normalization). Again this may be trivialized by fixing a defect line going through all $f_i$'s and modifying couplings along this line. Remark that $\prod_i f_i=1$ (otherwise the state space is empty).
For simplicity let us consider a pair $\chi(\sigma(v_2)),\chi^{-1}(\sigma(v_1)$ of spin variables, and a pair $\mu_g(f_2),\mu_{g^{-1}}(f_1)$ of disorder variables. Let us fix non intersecting paths $\gamma$ (resp. $\gamma^\dg$) on $\Gamma$ (resp. $\Gamma^\dg$) from $v_1$ to $v_2$ (resp. from $f_1$ to $f_2$). We can write a mixed correlator:
\begin{eqnarray*}
\lefteqn{{\mc Z}\langle \chi(\sigma(v_2))\chi^{-1}(\sigma(v_1))\mu_{g}(f_2)\mu_{g^{-1}}(f_1)\rangle=}
\\
& \hspace{1in} & \sum_{\sigma:V\rightarrow G}\prod_{(vv')\in \gamma}\chi(\sigma(v')\sigma(v)^{-1})\prod_{(vv')\in E}w'_e(\sigma(v')\sigma(v)^{-1})
\end{eqnarray*}
where $w'_e=w_e$ if $e$ does not cross $\gamma^\dg$ and $w'_e(.)=w_e(g.)$ otherwise (in which case we orient $e$ from right to left of $\gamma^\dg$). Repeating the argument above leads to
$$\langle  \chi_0(\sigma(v_2))\chi_0^{-1}(\sigma(v_1))\mu_{g}(f_2)\mu_{g^{-1}}(f_1)\rangle_\Gamma=\langle \mu_{\chi_0}(v_2)\mu_{\chi_0^{-1}}(v_1)\chi_{f_2}(g)\chi_{f_1}(g^{-1})\rangle_{\Gamma^\dg}$$
where the $\chi_f$'s are the dual order variables; weights on $\Gamma^\dg$ are Fourier coefficients of those on $\Gamma$. This extends to any number of insertions.

By analogy with the Ising model, on the square lattice we are particularly interested in self-dual weights, viz. weights such that $\widehat w\propto w\circ\phi$ for some isomorphism $\phi:\hat G\rightarrow G$. If we assume weights are nonnegative, by Parseval we have $\widehat w=w\circ\phi$.\\

{\bf Cyclic models.} \\

Let us specialize to the cyclic case: $G=\Z/q\Z$. $\U_q=\{z\in\C:z^q=1\}$. Let $\xi_0=\exp(2i\pi/q)$. Identifying $G$ with $\hat G$ in the usual way, the Fourier transform is written:
$$({\mc F}f)(j)=\frac 1{\sqrt q}\sum_{i=0}^{q-1}f(i)\xi_0^{-ij}$$
Finding a self-dual weight consists in finding eigenvectors for ${\mc F}$. This operator is the discrete Fourier transform (DFT), involved in fast Fourier transform algorithm. The question of its diagonalization is classical, if rather delicate. By Fourier inversion, the Fourier transform ${\mc F}:\C^G\rightarrow\C^G$ is of order 4: ${\mc F}^4=\Id$, and consequently the eigenvalues of ${\mc F}$ are $\{\pm 1,\pm i\}$. It is known that the multiplicity of the eigenvalue 1 is $\lfloor\frac q4\rfloor+1$.

First we look for self-dual Potts model weights. In this case: $w\propto 1+a\delta_0$. Plainly ${\mc F}w=\frac{a}{\sqrt q}+\sqrt q\delta_0$. 
The fixed point equation ${\mc Fw}=w$ is thus solved for $a=\sqrt q$. For Ising ($q=2$), we recover the self-dual weight function (up to multiplicative constant).

As soon as $q\geq 4$, the multiplicity of $1$ is $\geq 2$. Similarly to \cite{GrunbaumDFT}, one can look for a nice operator $D\in\End(C^G)$ (almost) commuting with ${\mc F}$. An educated guess gives:
$$D=(a+b\chi^{-1})R+(c+d\chi^{-1})$$
where $R$ is the right shift: $(Rf)(k)=f(k+1)$, and $\chi\in\C^G$ is the character: $\chi(k)=\xi_0^k$. Classically, ${\mc F}R=\chi{\mc F}$ ($\chi$ acting by pointwise multiplication) and ${\mc F}\chi^{-1}=R{\mc F}$. Besides $R\chi=\xi_0\chi R$. Then
$${\mc F}D=(a\chi+b\xi_0\chi R+c+dR){\mc F}$$
Consequently, ${\mc F}D=\lambda\chi D{\mc F}$ ($\lambda\in\C$) iff:
$$(a,b,c,d)=\lambda(b\xi_0,d,a,c)$$
that is: $\lambda^4=\xi_0^{-1}$, $(a,b,c,d)\propto(1,\lambda^3,\lambda,\lambda^2)$. An element $w$ of $\Ker D$ satisfies a trivial recursion: 
$$w(k+1)(1+\lambda^3\xi_0^{-k})+w(k)(\lambda+\lambda^2\xi_0^{-k})=0$$ 
To ensure positivity, we choose $\lambda=-\exp(-i\pi/2q)$; then we find
$$w(k)=\prod_{j=0}^{k-1}\frac{\sin(\frac{\pi k}{q}+\frac{\pi}{4q})}{\sin(\frac{\pi (k+1)}{q}-\frac{\pi}{4q})}$$
and check that $w(q)=w(0)=1$. Then $\Ker D$ is a line which is invariant under ${\mc F}$, and thus $w={\mc F}w$ is an invariant weight. This is the {\em Fateev-Zamolodchikov point} (\cite{FZZN}) for the $\Z/q\Z$ model (identified in \cite{FZZN} as part of a family of weights satisfying the Yang-Baxter relations).\\

{\bf Ashkin-Teller model.}\\

In the Ashkin-Teller model (\cite{Baxterexact}, as phrased by Fan), each vertex $v\in V$ carries a pair of spins $(\sigma_v,\rho_v)\in\{\pm 1\}^2$; the interaction weight on the edge $e=(vv')$ is:
$$\exp(\beta(J_e\sigma_v\sigma_{v'}+J_e'\rho_v\rho_{v'}+J''_e\sigma_v\sigma_{v'}\rho_v\rho_{v'}))$$
We then recognize the abelian spin model with $G=\{\pm 1\}^2$ (the simplest non cyclic finite abelian group). Identifying $G$ with $\hat G$ we get (here $(\eps_1,\eps_2)\in G$)
$$({\mc F} f)(\eps_1,\eps_2)=\frac 12\left(f(1,1)+\eps_1 f(-1,1)+\eps_2 f(1,-1)+\eps_1\eps_2 f(-1,-1)\right)$$
leading to the following simple self-duality condition for the weight $w:G\rightarrow\R$ (i.e. satisfying $w={\mc F}w$):
$$w(1,1)=w(-1,1)+w(-1,1)+w(-1,-1)$$
The basic order variables are $\sigma,\rho,\sigma\rho$ and there are naturally corresponding disorder variables $\mu,\nu,\mu\nu$. Remark that in the {\em isotropic} case $J=J'$, the model becomes a $\Z/4\Z$ model.\\

{\bf Parafermions.} \\

Let us go back to the general case of an abelian model with spin group $G$. Order (resp. disorder) variables are indexed by $\hat G$ (resp. $G$). By analogy with the Ising model, one can consider a {\em parafermion} (or {\em spinor})
$$\psi^{\chi_0,g_0}_{vf}=\chi_0(\sigma_v)\mu_{g_0}(f)$$
where $(\chi_0,g_0)\in\hat G\times G$ (and are fixed and omitted from now on), $v\in\Gamma$, $f\in\Gamma^\dg$ adjacent. When tracking an order variable $\chi_0(\sigma)$ along a cycle around a disorder variable $\mu_{g_0}(f)$ (or vice versa), one picks a phase $\chi_0(g_0)^{\pm 1}$. In particular, when the spinor rotates unto itself (eg fix $f$ and follow $v$ along $\partial f$), it is multiplied by $\chi_0(g_0)^{\pm 1}$. This may be formalized as follows. Consider a graph $\Gamma'$  with vertices at midpoints of segments $[vf]$, $v\in\Gamma$, $f\in\Gamma^\dg$, $v$ on the boundary of $f$ (denote such a vertex $(vf)$); edges of $\Gamma'$ are such that $(vf)\sim (v'f)$ if $v,v'$ are consecutive vertices on $\partial f$ and dually $(vf)\sim (vf')$ if $v\in\partial f\cap\partial f'$. The ``line bundle" $\bigoplus_{(vf)\in\Gamma'}\C(vf)$ is equipped with a connection (in the terminology of \cite{Kendouble} eg), i.e. a collection of isomorphisms $\phi_{u,u'}:\C u\rightarrow\C u'$ for adjacent vertices $u,u'$ in $\Gamma'$. These are defined by:
$$\phi_{(vf),(v'f)}(z)=\exp\left(is\arg\frac{f'-v}{f-v}\right)z,{\rm\ \ \ } \phi_{(vf),(vf')}(z)=\exp\left(is\arg\frac{v'-f}{v-f}\right)z$$
where $s$ is the spin, and arguments are chosen in $(-\pi,\pi)$ (so that $\phi_{u,u'}=\phi_{u',u}^{-1}$).

Let us consider a general correlator:
$$F(v_0,f_0)=\langle \psi_{v_0f_0}\prod_{i=1}^n\chi_i(v_i)\prod_{j=1}^m\mu_{g_j}(f_j)\rangle$$
We regard all insertions as fixed except $v,f$, which are adjacent. 
In order for this correlator to be non-trivial, we need $\chi_0\prod_{i\geq 1}\chi_i=1_{\hat G}$, $g_0\prod_{j\geq 1}g_j=1_G$. To assign an actual value to the correlator, we also have to fix a defect line (or defect tree) connecting all disorder operators on $\Gamma^\dg$. 
Alternatively, one can think of $F$ as being multivalued (with a phase ambiguity) or as a section of the line bundle described above (with the necessary modifications around $v_i,f_j$).

Consider the local situation around an edge $e=(vv')\in E$ and its dual edge $(ff')=(vv')^\dg\in E^\dg$; this gives four possible locations for the parafermion. We assume that the defect line $\gamma$ chosen for $f$ is the concatenation of $(ff')$ and the defect line $\gamma'$ for $f'$ (which has $f'$ as an endpoint).

For $g\in G$, let us consider the partial (modified) partition function:
$${\mc Z}_g=\sum_{\sigma:V\rightarrow G}\ind_{\sigma(v)=1,\sigma(v')=g}\prod_{i\geq 1}\chi_i(\sigma(v_i))\prod_{e=(xy)\in E, e\neq(vv')}w_e'(yx^{-1})
$$
where $w'$ is the weight function modified along $\gamma'$ in order to encode the disorder ${\mu_{g_0}(f')\prod_{j\geq 1}\mu_{g_j}(f_j)}$. Then up to a factor $|G|$
\begin{eqnarray*}
F(v,f)=\sum_{g\in G}w_e(gg_0){\mc Z}_g&&
F(v,f')=\sum_{g\in G}w_e(g){\mc Z}_g\\
F(v',f)=\sum_{g\in G}\chi_0(g)w_e(gg_0){\mc Z}_g&&
F(v',f')=\sum_{g\in G}\chi_0(g)w_e(g){\mc Z}_g
\end{eqnarray*}
Let us assume that the weight $w=w_e$ satisfies the following equation on $\C^G$:
\begin{equation}\label{eq:diff}
\left((a-b\chi_0)\theta_0-(c-d\chi_0)\right)w=0
\end{equation}
where $\theta_0$ is the shift: $(\theta_0 w)(g)=w(gg_0)$. Then
\begin{equation}\label{eq:para}
-a F(v,f)+b F(v',f)+c F(v,f')- d F(v',f')=0,
\end{equation}
the local {\em parafermionic equation} (\cite{SmiICM,Smiising,CarRiv,IkhCar}). Equation \ref{eq:diff} implies 
$$w(gg_0^k)=w(g)\prod_{j=1}^k\frac{c-d\chi_0(g)z^j}{a-b\chi_0(g)z^j}$$
where $z=\chi_0(g_0)$. In particular, if $r$ is the order of $g_0$
$$\prod_{n=0}^{r-1}(a-b z^n)=\prod_{n=0}^{r-1}(c-d z^n)$$
If $r'$ is the order of $z=\chi_0(g_0)$ in $\U$, we have: $\prod_{j=0}^{r'-1}(x-z^j)=x^{r'}-1$. Consequently we need:
$$\left(a^{r'}-b^{r'}\right)^{r/r'}=\left(c^{r'}-d^{r'}\right)^{r/r'}$$
Symmetry of the weight $w(g^{-1})=w(g)$ yields further conditions, such as:
$$\prod_{j=0}^{k}\frac{c-d z^j}{a-b z^j}\cdot\frac{c-d z^{-j-1}}{a-b z^{-j-1}}=1$$
coming from $w(g_0^k)=w(g_0^{-k})$ (assuming weights are nonzero). Hence
$$\frac{c-dz^k}{a-bz^k}=\frac{bz^{-1}-az^k}{dz^{-1}-cz^k}$$
Two homographies that agree on three points are equal. Thus if $r'\geq 3$, we have $(a,b,c,d)\propto (dz^{-1},c,bz^{-1},a)$, ie
$(a,b,c,d)\propto (u\lambda^{-1},u^{-1}\lambda,u^{-1}\lambda^{-1},u\lambda)$ where $u\in\C^*$, $\lambda^4=z$. In order for weights to be positive, we take $u$ unitary and recover the FZ weights (\cite{FZZN}):
$$w(g_0^k)=\prod_{j=0}^{k-1}\frac{\sin(\frac{\pi j}{r'}+\frac{\theta}{r'})}{\sin(\frac{\pi (j+1)}{r'}-\frac{\theta}{r'})}$$
with an additional ``anisotropy" parameter $\theta$ compared to the self-dual case. (In the cyclic case, arguing as before we find that ${\mc F}w_{\theta}\propto w_{\frac\pi 2-\theta}$).

\section{Discrete Gaussian Free Field}

We have considered abelian spin models with values in a finite abelian group. More generally one may consider a locally compact abelian group (the natural set-up for Fourier-Pontryagin duality). Let us consider scalar fields, i.e. fields taking values in $\R$.

As before, we consider a finite connected planar graph $\Gamma=(V,E)$. It is somewhat convenient to designate a non-empty subset $\partial$ of $V$ as the boundary. We consider discrete scalar fields, i.e. functions $\phi:V\rightarrow\R$. The values of $\phi$ on $\partial$ are fixed to some prescribed $\phi_{|\delta}$. This defines an affine state space. Let us also fix positive {\em conductances} $(c_e)_{e\in E}$.
 
The Discrete Gaussian Free Field (DGFF) on $(\Gamma,(c_e))$ with Dirichlet boundary condition $\phi_{|\delta}$ on $\partial$ is the Gaussian variable with measure proportional to
$$\prod_{e=(vv')\in E}\exp(-\frac {c_e}2(\phi(v')-\phi(v))^2)\prod_{v\in V\setminus\partial}\frac{d\phi(v)}{\sqrt{2\pi}}$$
For $f:V\rightarrow\R$, its Dirichlet energy can be written:
$$\frac 12\sum_{e=(vv')}c_e(f(v')-f(v))^2=\frac 12\sum_{v\in V} f(v)(\Lap_\Gamma f)(v)$$
by a discrete Green's formula argument, where $\Lap_\Gamma$ is the weighted graph (positive) Laplacian:
$$(\Lap_\Gamma f)(v)=\sum_{v'\sim v}c_{(vv')}(f(v)-f(v'))$$
Under our assumptions, 
$$\Lap_\Gamma:\{f:V\rightarrow\R,f_{|\partial}=0\}\rightarrow\{f:V\rightarrow\R,f_{|\partial}=0\}$$
is indeed invertible and the Dirichlet energy is a positive definite form on $\{f\in\R^V,f_{|\partial}=0\}$. Consequently the above Gaussian measure is finite. The mean field $\phi_0$ is the one minimizing the Dirichlet energy in the affine space $\{f\in\R^V,f_{|\partial}=\phi_{\partial}\}$, i.e. the solution of the Dirichlet problem: $(\Lap_\Gamma f)_{V\setminus\partial}=0$, $f_{|\partial}=\phi_{\partial}$.

The covariance kernel $C(v,v')=\Cov(\phi(v)\phi(v'))$ is the Green kernel $G$ for $\Lap_\Gamma$, viz. is characterized by $\Lap_\Gamma G(.,v')=\delta_{v'}$ in $V\setminus\partial$, $G(v,v')=0$ for $v\in\partial$. We deduce the characteristic function:
$$\E(\exp(i\sum_j\alpha_j\phi(v_j)))=\exp(i\sum_j\alpha_j\phi_0(v_j))\exp(-\frac 12\sum_{j,k}\alpha_j\alpha_k G(v_j,v_k))$$
The local variables $\exp(i\alpha\phi(v))$ are called {\em electric operators}, where $\alpha$ is the charge; these are the main order variables.

In order to introduce disorder variables, it is convenient to rephrase the problem in terms of the {\em current} $J(vv')=\phi(v')-\phi(v)$, a graph 1-form (in $\Omega^1(\Gamma)$). Let us also assume that there is no boundary. Then the Dirichlet energy is well-defined in terms of $J$: $\frac 12\sum_e c_e(J(e))^2$. As a state space, we consider the space of 1-forms such that $dJ=m$ where $m=\sum_j m_j\ind_{f_j}\in\R^{V^\dg}$ is a fixed magnetic charge distribution ($m_j$ is the {\em magnetic charge} positioned at the face $f_j$); here $dJ(f)=\sum_{\vec{e}\in\partial f}J({\vec e})$. The Gaussian variable $J$ induced on this affine space by the Dirichlet energy can be written $J=J_h+J_0$ where $J_h$, the mean, is the current of minimal energy with given magnetic charge distribution and $J_0$ is a centered field with covariance kernel:
$$\Cov(J(v_1v_2),J(v_3v_4))=G(v_2,v_4)+G(v_1,v_3)-G(v_2,v_3)-G(v_1,v_4)$$
(In the absence of a boundary, $G(.,v)$ is well defined only modulo an additive constant, which is enough to define the RHS).

If $J=d\phi$ is an exact 1-form, we have
$$\sum_{e\in E}c_eJ_e\omega_e=\sum_{v\in V}\phi(v)\sum_{v'\sim v}c_{vv'}\omega(vv')$$
which shows that the orthorgonal in $L^2(E,\sum c_e\delta_e)$ of exact 1-form are harmonic 1-forms, i.e. those satisfying $\sum_{v'\sim v}c_{vv'}\omega(vv')=0$ for all $v$. It is convenient to define a discrete Hodge star operator $\ast:\Omega^1(\Gamma)\rightarrow\Omega^1(\Gamma^\dg)$ defined by $(\ast\omega)((vv')^\dg)=c_{vv'}\omega(vv')$ (\cite{Mercat}). We regard $\Gamma^\dg$ as a weighted graph with $c_{e^\dg}=(c_e)^{-1}$, so that $\ast^2=-\Id$ and $\ast$ is isometric. Clearly, $\ast$ maps harmonic forms on $\Gamma$ to closed forms on $\Gamma^\dg$, and vice-versa. Thus $J\in\Omega^1(E)$ can be decomposed uniquely as $J=d\phi+\ast d\psi$ where $d\phi\in\Omega^1(\Gamma)$  and $d\psi\in\Omega^1(\Gamma^\dg)$ closed; this {\em Hodge decomposition} is orthogonal in $L^2(E,\sum c_e\delta_e)$. By analogy with the continuous case, we write $\langle J,J'\rangle_{L^2(E)}=\int J\wedge\ast J'$ (see \cite{Mercat}). Moreover, $\Lap_\Gamma=d\ast d$ (as usual identifying vertices of $\Gamma$ with faces of $\Gamma^\dg$). Consequently, the form $J_h$ with minimal Dirichlet energy among those such that $dJ=m$ is:
$$J_m=\ast d\left(\sum_{f\in V^\dg} m_fG^\dg(.,m_f)\right)$$
where $G^\dg$ is the Green kernel on $\Gamma^\dg$.\\

{\bf Duality.}\\

Let us consider abelian duality for the centered DGFF. At this point we could write mixed electric-magnetic correlators in terms of $G,G^\dg$ and see duality appear in these explicit formulae. For ease of comparison with the finite case, we phrase duality in terms of Fourier transform.

In general one can fix an electric charge distribution $(\alpha_v)_{v\in V}$ and a magnetic charge distribution $(m_f)_{f\in F}$, such that $\sum_v\alpha_v=0$, $\sum_f m_f=0$. For simplicity, let us consider two pairs of charges: electric charges $\pm\alpha$ at $v_{\pm}$, and magnetic charges $\pm m$ at $f_{\pm}$. If $v_{\pm}$ circles counterclockwise around $f_{+}$, the correlator picks up a factor $\exp(\pm i\alpha m)$; assuming $\alpha m\in 2\pi\Z$ removes multivaluedness issues. In the general case, let us fix paths $\gamma,\gamma^\dg$ on $\Gamma,\Gamma^\dg$ with endpoints $v_\pm, f_\pm$. We are interested in $\langle\exp(i\alpha\int_\gamma J)\rangle$, where $\langle.\rangle$ is the unnormalized expectation for the free field current $J$ on $\Gamma$ with mean $J_m$. We have the representation:
$$\int_\gamma d\phi=\int d\phi\wedge\ast dJ^\dg_\alpha$$
where $J^\dg_\alpha=\sum_{v\in V}\alpha_vG(.,v)$.

It is technically convenient to temporarily mollify the hard constraint $dJ=m$ and dampen the electric correlator by considering the Gaussian integral:
$$\int_{\R^E}\exp(i\alpha\int_\gamma d\phi-\frac\eps 2||d\phi||^2)\exp(-\frac{\eps^{-1}}2||\ast d\psi-J_m||^2)
\prod_{e\in E}\frac{e^{-\frac {c_e}2(J_e)^2}}{\sqrt {2\pi c_e^{-1}}}dJ_e$$
where $\eps$ is a small positive parameter, $J=d\phi+\ast d\psi$ ($\|.\|$ is the weighted $L^2$ norm, for which the Hodge decomposition is orthogonal). By Fourier inversion, this is written as
$$\int_{\R^E}\exp(-\frac {\eps^{-1}}2||d\tilde\phi-dJ^\dg_\alpha||^2)\exp(i\int d\tilde\psi\wedge dJ_m-\frac\eps 2||d\tilde\psi||^2)\prod_{e\in E}\frac{e^{-\frac {c_e^{-1}}2(\tilde J_e)^2}}{\sqrt {2\pi c_e}}d\tilde J_e$$
where $\tilde J=d\tilde\phi+\ast d\tilde\psi$ is the Hodge decomposition on $\Gamma^\dg$.

Letting $\eps\searrow 0$, we conclude that the discrete free fields on $\Gamma,\Gamma^\dg$ (with reciprocal conductances) are in duality, in which magnetic and electric charges are exchanged.\\

{\bf Compactification.}\\

In the case where the underlying graph $\Gamma$ is embedded on a torus $\Sigma=\C/L$, if we consider the Gaussian measure on closed currents $J: dJ=0$ induced by the Dirichlet energy, we gain two ``topological" marginals: $(\int_A J,\int_B J)$, where $A,B$ are two standard cycles generating the homology of $\Sigma$. Given $\alpha,\beta$, there is a unique closed, harmonic form $d\omega_{\alpha,\beta}$ ($d\omega_{\alpha_\beta}=0$, $d\ast \omega_{\alpha,\beta}=0$) such that $(\int_A \omega_{\alpha,\beta},\int_B \omega_{\alpha,\beta})=(\alpha,\beta)$. The space of 1-forms decomposes as an orthogonal sum of exact, coexact, and closed and harmonic forms: this is (a discrete version of) the Hodge decomposition. Thus an instance of a free field current $J$ can be decomposed as $J=d\phi+\omega$, where $\phi$ is a scalar field (on $V$, defined modulo additive constant) and $\omega$, the instanton component, is a closed harmonic form; these two summands are independent. 

It is classical (\cite{GawCFT,DiF}) to consider a {\em compactified} version of the field, in which in the periods $(\int_A J,\int_B J)$ are constrained to take values in a prescribed lattice $2\pi r\Z$ ($r$ is the {\em compactification} radius). This has no effect on the scalar component of the field; the instanton becomes a discrete variable supported on $\{\omega_{2\pi mr,2\pi nr}\}$, with weights proportional to the Dirichlet energy of these forms (a ``discrete Gaussian" variable).

If $\Lambda$ is a lattice in a Euclidean space $(V,\langle.,.\rangle)$, its dual lattice $\Lambda'$ ($y\in\Lambda'$ iff $\langle x,y\rangle\in\Z$ for all $x\in\Lambda$), we have the following version of the Poisson summation formula:
$$\sum_{x\in\Lambda}f(x)={\rm Vol}(V/\Lambda)^{-1}\sum_{y\in\Lambda'}\hat f(y)$$
in the normalization $\hat f(y)=\int_Ve^{-2i\pi\langle y,x\rangle} f(x)dx$. In particular for $f(x)=e^{-\pi t\|x\|^2}$, we get
$$\sum_{x\in\Lambda}e^{-\pi t\|x\|^2}=t^{-\dim V/2}{\rm Vol}(V/\Lambda)^{-1}\sum_{y\in\Lambda'}e^{-\pi t^{-1}\|y\|^2}.$$
Let us specialize this to $V$ the two-dimensional space of closed harmonic forms (which represents $H^1(\Sigma)$) with norm given by $2\pi$ times the Dirichlet energy; $\Lambda=\{\omega:\int_{A,B}\omega\in 2\pi r\Z\}$. In order to identify $\Lambda'$, we introduce the closed harmonic 1-forms $\omega^\dg_{c,d}$ on $\Gamma^\dg$ such that $\int_{A,B}\omega^\dg_{c,d}=c,d$ and recall the bilinear relation:
$$\langle \omega_{a,b},\ast\omega^\dg_{c,d}\rangle=\int_\Sigma \omega_{a,b}\wedge\omega^\dg_{c,d}=\oint_{\partial F}h_{a,b}\omega^{\dg}_{c,d}=ad-bc$$
where $F$ is a fundamental domain bounded by the cycles $a,b$ and $dh_{a,b}=\omega_{a,b}$ in $F$. We conclude $\ast\Lambda'=\{\omega^\dg:\int_{A,B}\omega^\dg\in 2\pi r^{-1}\Z\}$. Thus if we write the {\em instanton partition function}
$${\mc Z}_{inst}(r)=\sum_{\omega\in\Lambda}e^{-\frac 12\int \omega\wedge\ast\omega}$$
which represents the (multiplicative, by independence) contribution of the instanton component to the partition of the compactified free field at radius $r$, we have 
$${\mc Z}_{inst}(r)={\mc Z}_{inst}^\dg(r^{-1})$$
up to an elementary multiplicative factor. Thus duality on compactified free fields has the effect of inverting the compactification radius (\cite{GawCFT}). This is a simple example of {\em $T$-duality}.\\

{\bf Gaussian free field.}\\

The continuous Gaussian (or massless) free field on, say, a torus $\Sigma$ is the Gaussian field with action
$$S(\phi)=\frac{g}{4\pi}\int_\Sigma |\nabla \phi|^2$$
($g>0$ the {\em coupling constant}), i.e. the centered Gaussian field with covariance kernel $\frac{2\pi}{g}(-\Lap)^{-1}$ where we consider $(-\Lap)$ as an invertible operator on zero-mean functions on $\Sigma$ (\cite{SimonPphi,GlimmJaffe}). In the abstract Wiener space approach, one thinks of the free field as a random element of a large enough Banach space, typically a negative index Sobolev space. Plainly it may also be considered on the plane, finite domains, Euclidean spaces of other dimensions and Riemannian manifold. The compactified free field with values in $\R/2\pi r\Z$ is defined as above as the sum of a scalar free field and an instanton component, which is supported on a lattice of harmonic 1-forms.

Let us briefly discuss magnetic and electric charges in the plane (\cite{DiFSalZubcoulomb},\cite{GawCFT}). In this case, the covariance kernel is $C(z,w)=g^{-1}\log|z-w|$ (with a slight abuse of terminology, as we have to quotient by constant functions). The electric operators are written formally ${\mc O}_e(z)=\exp(ie\phi(z))$ ($e$ is the electric charge); this is somewhat problematic, as a realization of the field $\phi$ is only a distribution. In the standard regularization, one disregards the diverging self-energy and write:
$$\langle:{\mc O}_e(z){\mc O}_{-e}(w):\rangle=|z-w|^{-\frac{e^2}g}$$
where the colons denote the regularization procedure. Alternatively, one may replace $\phi(z)$ with its average on $D(z,\delta)$ and define
$$\langle:{\mc O}_e(z){\mc O}_{-e}(w):\rangle=\lim_{\delta\rightarrow 0}\delta^{-\frac{e^2}g}\left\langle\exp(ie\left(\frac{1}{\pi\delta^2}\int_{D(z,\delta)}\phi-\frac{1}{\pi\delta^2}\int_{D(w,\delta)}\phi\right)\right\rangle$$
In the presence of a magnetic charge $m$ at $z$ (represented by ${\mc O}_m(z)$), the field increases by $2m\pi$ when tracked along a counterclockwise around $z$. For a pair of insertions ${\mc O}_m(z){\mc O}_{-m}(w)$, the average field is $m(\arg(.-z)-\arg(.-w))$, the harmonic conjugate of $m(\log|.-z|-\log|.-w|)$. Its Dirichlet energy has a logarithmic blow-up at singularities; a standard finite part regularization gives the expression
$$\langle:{\mc O}_e(z){\mc O}_{-e}(w):\rangle=|z-w|^{-m^2g}$$
Combining these two elements, we obtain:
\begin{eqnarray*}
\lefteqn{
\langle:{\mc O}_e(z_2){\mc O}_{-e}(z_1){\mc O}_m(w_2){\mc O}_{-m}(w_1):\rangle=}\\
& & |z_1-z_2|^{-\frac{e^2}g}|w_1-w_2|^{-m^2g}e^{ime(\arg(z_2-w_2)-\arg(z_2-w_1)-\arg(z_1-w_2)+\arg(z_1-w_1))}
\end{eqnarray*}
To obtain ``spinor" variables ${\mc O}_{em}$, we coalesce electric and magnetic charges: ``${\mc O}_{em}(z)=\lim_{z\rightarrow w}{\mc O}_e(z){\mc O}_m(w)$". Specifically take $z_i=w_i+\delta u_i$, $|u_i|=1$, $\delta\searrow 0$. Then
$$\langle:{\mc O}_{em}(w_2){\mc O}_{-e,-m}(w_1):\rangle=|w_2-w_1|^{-\frac{e^2}g-m^2g}e^{-2iem\arg(w_2-w_1)}e^{iem(\arg(u_1)+\arg(u_2)}$$
The last part is a manifestation of the spinor nature of these variables. Fixing a reference direction for $u_1,u_2$, this may be rewritten as:
$$\langle:{\mc O}_{em}(w_2){\mc O}_{-e,-m}(w_1):\rangle=(w_2-w_1)^{-\frac{e^2}{2g}-\frac{m^2g}2-em}(\overline{w_2-w_1})^{-\frac{e^2}{2g}-\frac{m^2g}2+em}$$

In this context, one recovers the duality between (properly renormalized) electric and magnetic operators ${\mc O}_{em}\leftrightarrow{\mc O}_{2m,e/2}$ (and $g\leftrightarrow \frac 4g$ for couplings), and pairs of reciprocal radii in the compactified case (\cite{GawCFT,DiF}).

\section{Random-cluster model}

Let us consider the Ising model on a graph $\Gamma=(V,E)$. The weight associated to the edge $e$ and the spin configuration $(\sigma_v)_{v\in V}$ is $\exp(\beta J_e\sigma_{v}\sigma_{v'})$. As we have seen, writing this weight as
$$\exp(\beta J_e\sigma_{v}\sigma_{v'})=\cosh(\beta J_e)\left(1+\sigma_{v}\sigma_{v'}\tanh(\beta J_e)\right)
$$
leads to the high-temperature graphical representation of the Ising model. It may also be written as
$$\exp(\beta J_e\sigma_{v}\sigma_{v'})=e^{-\beta J_e}\left(1+\delta_{\sigma_v,\sigma_{v'}} (e^{2\beta J_e}-1)\right)
$$
which leads to the {\em random-cluster representation} (\cite{GrimFK}), as we now explain. Set $w(e)=e^{2\beta J_e}-1$. The partition function reads:
\begin{align*}
{\mc Z}&=\sum_{\sigma:V\rightarrow\{\pm 1\}}\exp(\beta J_e\sigma_{v}\sigma_{v'})=\prod_{e\in V}e^{-\beta J_e}\sum_{\sigma:V\rightarrow\{\pm 1\}}\left(1+\delta_{\sigma_{v},\sigma_{v'}}w(e)\right)\\
&=\prod_{e\in V}e^{-\beta J_e}\sum_{\sigma:V\rightarrow\{\pm 1\}}\sum_{E_0\subset E}\prod_{e=(vv')\in E_0}\delta_{\sigma_{v},\sigma_{v'}}w(e)
\end{align*}
We notice that $\sum_{\sigma:V\rightarrow\{\pm 1\}}\prod_{(vv')\in E_0}\delta_{\sigma(v),\sigma(v')}=2^{C(E_0)}$, where $C(E_0)$ is the number of connected components of the subgraph $\Gamma_0=(V,E_0)$. Hence
$${\mc Z}=\prod e^{-\beta J_e}\sum_{E_0\subset E}2^{C(E_0)}\prod_{e\in E_0}w(e)$$
More generally, one may consider a $q$-state Potts model with edge weight $e^{2\beta J_e\delta_{\sigma_{v},\sigma_{v'}}}$, which leads to the following random cluster representation of the partition function:
$${\mc Z}=\sum_{E_0\subset E} q^{C(E_0)}\prod_{e\in E_0}w(e)$$
with $w(e)=e^{2\beta J_e}-1$. Plainly, the RHS may be interpreted as the partition function of a model where the configuration space consists of subgraphs $\Gamma_0=(V,E_0)$ of $\Gamma$ and the configuration weight of $\Gamma_0$ is proportional to
$$q^{C(E_0)}\prod_{e\in E_0}w(e)$$
This is the {\em random-cluster model} or {\em Fortuin-Kasteleyn (FK) percolation} (\cite{GrimFK}), and is defined for any $q>0$ (whereas the Potts model is restricted to $q\in\N^*$). 

One may also track natural correlators, such as spin correlators, through the random-cluster representation. We see that in the $q$-state Potts model,
$$\langle \sigma(v_2)\sigma(v_1)^{-1}\rangle_{{\rm Potts}}=\P_{FK}(v_1\leftrightarrow v_2)
$$
where $v_1\leftrightarrow v_2$ means that $v_1,v_2$ belong to the same connected component (cluster) in the random subgraph induced by FK-percolation, and Potts spins are regarded as elements of $\U_q=\{z\in \U:z^q=1\}$. By the ${\mf S}_q$ symmetry of the Potts model, we have $\langle\delta_{\sigma(v),1}\rangle_{{\rm Potts}}=\frac 1q$ and $\langle \sigma(v_2)\sigma(v_1)^{-1}\rangle_{{\rm Potts}}=p+\frac{1-p}{1-q}$ where $p=\langle \delta_{\sigma(v_1),\sigma(v_2)}\rangle$; then
$$\Cov(\delta_{\sigma(v_1),1},\delta_{\sigma(v_2),1})
=\frac pq-\frac 1{q^2}=\frac{q-1}{q^2}\langle \sigma(v_2)\sigma(v_1)^{-1}\rangle_{{\rm Potts}}=\frac{q-1}{q^2}\P_{FK}(v_1\leftrightarrow v_2)
$$

Higher order spin correlations are easiest to understand with the Edwards-Sokal coupling (\cite{EdwSok}) between the $q$-Potts model and the corresponding FK percolation. Starting from an FK configuration, one chooses a spin uniformly at random for each cluster and assign it to each vertex in this cluster. Then the resulting spin configuration is $q$-Potts distributed. Consequently, a general spin correlator $\langle\prod_j\chi_j(\sigma(v_j))\rangle_{{\rm Potts}}$ may be expressed in terms of the connectivity of $v_1,\dots,v_n$ in the FK percolation subgraph (which induces a random partition of $\{v_1,\dots,v_n\}$).\\

{\bf Duality.}\\

In order to discuss duality, we need to specialize to planar graphs. In percolation and related models (\cite{Grimperc}), it is standard to introduce a dual configuration $E_0^\dg\subset E^\dg$ defined by $e^\dg\in E_0^\dg$ iff $e\notin E_0$, so that $|E_0|+|E_0^\dg|=|E|$. Euler's formula for $\Gamma_0=(V,E_0)$ reads:
$$|V|-|E_0|+|F_0|=1+C(E_0)$$
where $F_0$ designates the faces of $\Gamma_0$ (on the Riemann sphere) and is easily seen to be in bijection with the connected components of $\Gamma_0^\dg=(V^\dg,E_0^\dg)$. Hence
$$q^{C(E_0)}\prod_{e\in E_0}w(e)=\left(q^{|V|-|E|-1}\prod_{e\in E}w(e)\right)q^{C(E_0^\dg)}\prod_{e^\dg\in E_0^\dg}\frac{q}{w(e)}$$
Hence by defining a dual edge weight $w(e^\dg)=\frac{q}{w(e)}$, we see that $E_0\leftrightarrow E_0^\dg$ is a measure preserving correspondence between FK percolation on $(\Gamma,(w(e))_{e\in E})$ and $(\Gamma,(w(e^\dg))_{e^\dg\in E^\dg})$. 

In the case $q\in\{2,3,\dots\}$, this is consistent with the abelian duality for the Potts model, seen as a $\Z/q\Z$ model. Indeed, if the edge weight function $w_e:G\rightarrow\R^+$ is:
$$w_e\propto 1+w_e\delta_1$$ 
its Fourier transform is:
$${\mc F}w_e\propto \frac{w_e}{\sqrt q}+\sqrt q\delta_1$$
yielding back $w(e^\dg)=\frac q{w(e)}$.\\

{\bf Loop representation.}\\

In turn, the random cluster representation (in the planar case) may be mapped to a {\em dense loop} (or fully-packed loop) representation (in the standard case of the square lattice with self-dual weights, this is the {\em dense $O(n)$ model} at $n=\sqrt q$), which we now describe. Attached to the pair $(\Gamma,\Gamma^\dg)$ of graphs on the Riemann sphere, one constructs the quadrangulation $\dmd$ (the ``diamond graph", \cite{Mercat}) with vertices $V_{\dmd}=V\sqcup V^\dg$ and edges $E_\dmd$ defined by $(vf)\in\dmd$ is $v$ is a vertex of $V$ on the boundary of $f\in V^\dg$ (identified with faces of $\Gamma$). Then $D\dmd$, the derived (or medial) graph of $\dmd$, is specified by $V_{D\dmd}=E_\dmd$ and $(e_1e_2)\in E_\dmd$ if $e_1,e_2$ are two consecutive edges on a face of $\dmd$. (This is a planar version of the derived graph). Note that exchanging the roles of $\Gamma,\Gamma^\dg$ leads to the same $\dmd$, $D\dmd$.

To a pair of dual configurations $E_0\subset E$, $E_0^\dg\subset E^\dg$, one associates $E'_0\subset E_{D\dmd}$, the set of all edges in $D\dmd$ that do not intersect edges in $E_0,E_0^\dg$. Specifically, if $(v_1v_2)\in E$, $(f_1f_2)\in E^\dg$, then either $(v_1v_2)\in E_0$ and thus $((f_1v_1),(f_1v_2))$, $((f_2v_1),(f_2v_2))\in E'_0$; or $(f_1f_2)\in E_0^\dg$ and thus $((f_1v_1),(f_2v_1))$, $((f_1v_2),(f_2v_2))\in E'_0$. By construction $E'_0$ is 2-regular, hence consists in a set of disjoint closed loops that cover $V_{D\dmd}$ (a {\em loop gas}; note that not all spanning collections of loops on $D\dmd$ are obtained in this way). These loops separate clusters in $\Gamma_0$ from clusters in $\Gamma_0^\dg$, and it is easily seen that the number of loops is $C(E'_0)=C(E_0)+C(E_0^\dg)-1$. 
As $V-|E_0|=C(E_0)-C(E_0^\dg)+1$, we have $2C(E_0)=|V|-|E_0|+C(E'_0)$. Hence the FK weight of the configuration $E_0$ may be written as:
$$q^{C(E_0)}\prod_{e\in E_0}w(e)\propto \sqrt{q}^{C(E'_0)}\prod_{e\in E_0}w(e)q^{-\frac 12}$$
Consequently, up to normalization, an FK configuration on $\Gamma$ (or the dual configuration on $\Gamma^\dg$) maps to a loop gas $E'_0$ on $D\dmd$ with weight
$$\sqrt q^{C(E'_0)}\prod_{e'\in E'_0}w'(e')$$
with $w'(e')=\left(\frac{w(vv')}{\sqrt q}\right)^{\frac 14}$ if $e'=((vf),(v'f))$; and $w'(e')=\left(\frac{w(ff')}{\sqrt q}\right)^{\frac 14}$ if $e'=((vf),(vf'))$ (recall that $w(ff')w(vv')=q$ if $(ff')=(vv')^\dg$).\\

{\bf Parafermions.}\\

In the $q$-Potts case, Nienhuis and Knops (\cite{NieKnopotts}) showed that parafermion correlators could be computed in terms of the loop gas configuration. Consider $v_i$ a vertex of $\Gamma$ on the the boundary of the face $f_i$, $i=1,2$. One may introduce disorder operators $\mu_{-k}(f_1)$, $\mu_{k}(f_2)$, with $k\in\Z/q\Z\neq 0$ fixed. The new configuration space consists of $\Z/q\Z$ valued 1-forms $\omega$ on $\Gamma$ with $d\omega=k\delta_{f_2}-k\delta_{f_1}$. One may apply the random cluster representation to this state space: if $e\in E$ is such that $\omega(e)=1$, then $e\in E_0$ with probability $p_e$, where $\frac{p_e}{1-p_e}=w(e)$. The resulting random clusters may {\em not} wrap around $f_1,f_2$. In this case, one also obtains a modified (less canonical) Edwards-Sokal coupling. For instance, one may consider a reduced graph the vertices of which are the random clusters, with edges corresponding to edges of $E$ s.t. $\omega(e)\neq 1$. Then the number of assignments compatible with the disorder condition $q^{|C|-1}$. One then assigns a spin uniformly at random on each cluster, and set a defect line from $f_1$ to $f_2$ which does not cross any cluster (i.e. is drawn on a dual cluster). There is an overall factor $q$ depending on whether one quotients the state space by the global $\Z/q\Z$ symmetry.

We now wish to evaluate $\langle \sigma(v_2)^\ell\sigma(v_1)^{-\ell}\mu_{k}(f_2)\mu_{-k}(f_1)\rangle$, where $\sigma(v)\in \U_q=\{1,e^{\frac{2i\pi }q},\dots, e^{\frac{2i\pi (q-1)}q}\}$ is the spin variable, and $\ell\in\Z/q\Z$, $k,\ell\neq 0$. Plainly, the only contributing RC configurations are those in which $v_1,v_2$ are connected by a cluster and $f_1,f_2$ are connected by a dual cluster. Equivalently, in the loop representation, there is a loop $L$ passing through $(v_1f_1)$ and $(v_2f_2)$. Configurations containing $L$ will be counted with a phase depending of the isotopy type of $L$ in $\C\setminus\{v_1,f_1,v_2,f_2\}$. We now explain how to compute this phase. For definiteness, let us fix a defect line $\gamma$ from $f_1$ to $f_2$. 

Consider the half-loop $L^+$ (resp. $L^-$) of $L$ starting from $(v_1f_1)$  and ending at $(v_2f_2)$ with $v_1,v_2$ on its LHS (resp. RHS). Let us orient both half-loops from $(v_1f_1)$ to $(v_2f_2)$. One may also orient all other loops arbitrarily. If all loops have an orientation, one may define a height function $h:V_\dmd\rightarrow\frac\pi 2\Z$ as follows: if $v,f$ are adjacent vertices in $V_\dmd$, $h(f)-h(v)=\pm\frac\pi 2$ according to whether the loop crossing $\overrightarrow{vf}$ crosses it from left to right or right to left; then the loops may be seen as level sets of $h$. In the present case, the definition is ambiguous due to the conflicting orientations of $L^+,L^-$. Then $h$ may be seen as an additively multivalued function, picking a constant $\pm\pi$ when circling around $(v_1f_1)$, $(v_2f_2)$; this is closely analogous to the notion of magnetic charge we discussed for the DGFF.

Take a path from $v_1$ to $v_2$ on the cluster that contains them. It crosses $\gamma$ $n$ times (algebraically: crossing $\gamma$ from left to right counts +1, from right to left counts -1); and $\sigma(v_2)\sigma(v_1)^{-1}=\exp(\frac{2i\pi kn}q)$. Specifically, we can look at the leftmost (resp. rightmost) such path, that tracks $L^+$ (resp. $L^-$) on its RHS (resp. LHS). Each time one of these crosses $\gamma$, it corresponds to $L^\pm$ crossing $\gamma$, which contributes to $\int_\gamma dh$ (the height variation $h(f_2)-h(f_1)$ evaluated along $\gamma$). Hence $\pi n=\int_\gamma dh$, and:
$$\langle \sigma(v_2)^\ell\sigma(v_1)^{-\ell}\mu_{k}(f_2)\mu_{-k}(f_1)\rangle_{{\rm Potts}}=u\left\langle \ind_{(v_1f_1)\leftrightarrow (v_2f_2)}\exp\left(\frac{2ik\ell} q \int_\gamma dh\right)\right\rangle_{{\rm O(n)}}$$
where $|u|=1$ can be made explicit in terms of $\gamma$, $\arg\overrightarrow{f_1v_1}$, $\arg\overrightarrow{f_2v_2}$. In the case where $k$ or $\ell$ generates $\Z/q\Z$, it is clear that the LHS depends on $k,\ell$ only through $k\ell$.

Notice that the expression
$$\left\langle \ind_{(v_1f_1)\leftrightarrow (v_2f_2)}\exp\left(2is\int_\gamma dh)\right)\right\rangle_{{\rm O(n)}}$$
is defined for any $s\in\R/\Z$, $q>0$, and may thus be used to define parafermionic correlators for general FK percolation (\cite{SmiICM,Smiising,CarRiv}).

\section{The 6-vertex model}

\begin{figure}[htb]
\begin{center}
\scalebox{0.5}{\includegraphics{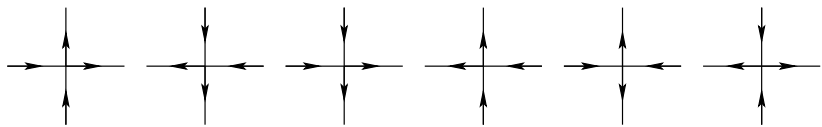}}
\end{center}
\caption{6V configurations}
\label{Fig:6Vconf}
\end{figure}

The 6-vertex (6V model for short) can be phrased for any finite 4-regular graph $\Gamma=(V,E)$, such as the square lattice (or some portion of it) or the Kagom\'e lattice. Its free energy per site was evaluated by Lieb, based on a Bethe ansatz analysis (\cite{Liebsquare,Baxterexact,Resh6V}). A configuration consists in an assignment of orientation (``arrow") to each edge, in such a way that every vertex is the origin and the endpoint of two oriented edges. Thus there are ${4\choose 2}=6$ possible local configurations around a given vertex (Figure \ref{Fig:6Vconf}). Graphs with higher (even) valencies may be considered: for a 6-regular graph (eg the triangular lattice), one may consider an orientation such that each vertex has 3 ingoing and 3 outgoing edges. This defines the (configuration space of) the 20-vertex model (as ${6\choose 3}=20$). Note that reverting all orientations preserves the 6V condition.

In the 6V model, one associates a Gibbs weight to each configuration, which is a product of local factors, read from the local configuration around each vertex. In the planar case (or oriented case), one may order edges cyclically around a vertex; then there are four (in-in-out-out) configurations (up to rotation: types 1-4) and two (in-out-in-out) configurations (types 5-6). Types are also numbered so that 1-2,3-4,5-6 are exchanged by reversal of all orientations. 

Let us denote $w_1(v),\dots,w_6(v)$ the local Gibbs weights around the vertex $v$ and $I_v=I_v(\omega)\in\{1,6\}$ the type at $v$ of the configuration $\omega$; then the Gibbs weight of $\omega$ is:
$$w(\omega)=\prod_{v\in V}w_{I_v(\omega)}(v)$$
Let us denote $\bar\omega$ the configuration obtained by reversing all arrows in $\omega$. One usually requests:
$$w_1=w_2=a,w_3=w_4=b,w_5=w_6=c$$
so that $w(\bar\omega)=w(\omega)$. On the square lattice, the rotationally invariant weights ($w_1=\cdots=w_4$, $w_5=w_6$) give the {\em $F$-model}.

Remark that replacing the weights $(w_i(v))$ with $(\lambda_vw_i(v))$ has the effect of multiplying weights of all configurations by a constant $\prod_{v\in V}\lambda_v$, a trivial modification. Thus there are essentially two free parameters for the weights at $v$: an anisotropy parameter $a/b$, and another parameter:
$$\Delta=\frac{a^2+b^2-c^2}{2ab}$$

We see the underlying graph $\Gamma$ as embedded on an oriented surface $\Sigma$. Then there is a dual (embedded) graph $\Gamma^\dg$, such that each vertex of $\Gamma^\dg$ corresponds to a face (assumed contractible) of $\Gamma$, and vice versa; edges of $\Gamma$ and $\Gamma^\dg$ are in bijection. Associated to a 6V configuration $\omega$ on $\Gamma$ is a (possibly additively multivalued) function $h=h(\omega)$ on $\Gamma^\dg$, defined as follows. If $f,f'$ are faces of $\Gamma$ (vertices of $\Gamma^\dg$) and $(vv')$ is the edge of $\Gamma$ between $f,f'$ {\em with orientation given by} $\omega$, then
\begin{align*}
h(f')-h(f)&=\frac\pi 2&{\rm\ if\ }\overrightarrow{vv'}{\rm\ crosses\ }\overrightarrow{ff'}{\rm\ from\ left\ to\ right}\\
h(f')-h(f)&=-\frac\pi 2&{\rm\ if\ }\overrightarrow{vv'}{\rm\ crosses\ }\overrightarrow{ff'}{\rm\ from\ right\ to\ left}
\end{align*}
The 6V rule ensures that $\sum_{i=0}^3 h(f_{i+1})-h(f_{i})=0$ if $f_0,\dots,f_3$ are the four faces around $v\in V$, cyclically indexed. This identifies the 6V model as a body-centered solid-on-solid (BCSOS) model (see \cite{VB6VSOS}). If $\Sigma$ is contractible, $h$ is well defined, up to an additive constant. Otherwise, for instance when $\Sigma$ is a torus, $h$ is additively multivalued, i.e. picks an additive constant (an integer multiple of $\pi$) when traced along a non-contractible cycle on $\Sigma$. More precisely, the graph 1-form $J=dh:E^\dg\rightarrow\R$ given by
$$J((ff'))=h(f')-h(f)=-J((f'f))$$
is well-defined and closed:
$$dJ(v)=\sum_{e\in\partial v}J(e)=0$$
where $v$ is a vertex of $\Gamma$ (a face of $\Gamma^\dg$) and $\partial v$ are the edges bounding $v$ in $\Gamma^\dg$, oriented counterclockwise. If $\gamma$ is a closed loop on $\Gamma^\dg$, $\int_\gamma J=\sum_{e\in\gamma}J(e)$ depends on $\gamma$ only through its homology class in $H_1(\Sigma)$. 

We may then define several natural order and disorder variables, by analogy with free field observables. First we have the current variables $J(e^\dg)\in\pm\frac\pi 2$ (which are sometimes thought of as $\pm 1$ spin variables). Then one can define natural electric correlators of type:
$$\langle\exp(i\alpha(h(f')-h(f)))\rangle$$
(for definiteness, one needs to fix a path from $f$ to $f'$), with $\pm\alpha$ electric charge at $f,f'$.

A convenient way to introduce disorder (magnetic) operators is to introduce edge defects. An edge $e$ split into two half-edges with outward orientation represents a magnetic charge $\frac 12$ (i.e. the height field picks up $\frac 12(2\pi)$ when cycling counterclockwise around this charge); symmetrically, an edge $e$ split into two half-edges with inward orientation represents a magnetic charge $-\frac 12$.

Let us now restrict to the standard case of the square lattice with periodic conditions ($\frac 1n\Z^2/L$, where $L$ is a sublattice, $L\sim \Z+\tau\Z$ where $n\gg 1$ is a scale parameter and $\tau\in\H=\{z:\Im(z)>0\}$ fixed), and isotropic weights, so that we may set $a=b=1$, and $\Delta=1-\frac{c^2}2$. We also assume that noncontractible cycles have even length. 

The main conjecture on the fluctuations of the six-vertex model, based on a Coulomb gas representation (\cite{NienCoulomb}), is that in the small mesh limit ($n\rightarrow\infty$), the height field converges to a compactified free field with compactification radius $\frac 12$ and coupling constant
$$g=\frac 8{\pi}\arcsin(\frac{c}2)$$
i.e. with action:
$$S(\phi)=\frac g{4\pi}\int_\Sigma|\nabla\phi|^2$$
where $\Sigma=\C/\tau\Z$ is the limiting torus. This is expected to hold in the ``regime III" phase, i.e. when $-1<\Delta<1$ or $0<c<2$.\\

{\bf Relation with FK percolation.}\\

For clarity we stick to the square lattice with periodic boundary conditions (and even periods). There is a close relation between FK percolation and the six-vertex model, which is best expressed in terms of the {\em Baxter measure}, which we now describe.

\begin{figure}[htb]
\begin{center}
\scalebox{.5}{\includegraphics{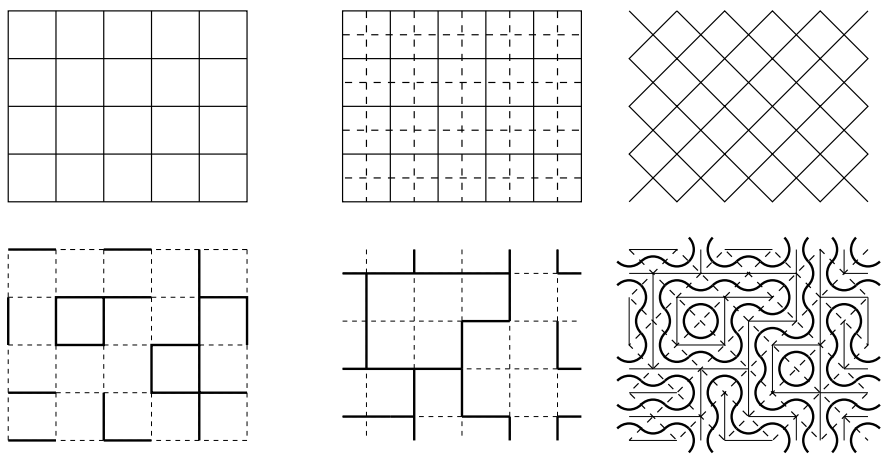}}
\end{center}
\caption{First row: graph $\Gamma$ (with periodic boundary conditions); its dual $\Gamma^\dg$ (dashed); $\dmd$. Second row: a primal FK configuration; dual configuration; dense $O(n)$ configuration
}
\label{Fig:FKON}
\end{figure}

Starting from $\Gamma$ modelled on the square lattice, we consider 
its dual graph $\Gamma^\dg$ and the diamond graph $\dmd$, which is itself locally a (scaled, rotated) square lattice. Associated to an FK configuration on $\Gamma$ and its dual configuration, we have constructed a dense $O(n)$ configuration (loop gas) on $D\dmd$. It is convenient to represent its segments as quarter circles: in each face of $\dmd$, two nonintersecting quarter circles connect the four edge midpoints. Starting from the self-dual FK weights, the weight of the $O(n)$ is simply (up to normalization) $n^\ell$ where $\ell$ is the number of loops and $n=\sqrt q$. See Figure \ref{Fig:FKON}.

The $O(n)$ state space may be enriched by looking at configurations of {\em oriented} loops. The winding of each loop if $\pm 2\pi$ for contractible loops (which are simple) and $0$ for non-contractible loops (on the torus). Each quarter circle is either a left or right turn and 
the weight of a oriented loop configuration is set to 
$$\exp(is\frac\pi 2 ({|{\rm right\ turns}|-|{\rm left\ turns}|}))$$
where $s$ is real. Collecting all weights along a given loop $L$ gives $e^{2i\pi s.wind(L)}\in\{\chi^{-1},1,\chi\}$, where $\chi=e^{2i\pi s}$. 

Remarkably, this complex measure on oriented loops has two real projections: one is a six-vertex model and the other is a weighted FK percolation measure. 

Let us start with the latter: by forgetting orientations, oriented loop configuration configurations map to dense $O(n)$ configurations, with weight per loop: $n=\chi+\chi^{-1}$ for contractible loops and $n=1$ for non-contractible loops. In the thermodynamic limit, one expects the number of these non-contractible loops (or, for that matter, the cardinality of any class of macroscopic loops) to stay tight. If $0<q<4$, one chooses $s$ so that $2\cos(2\pi s)=\sqrt{q}$; then the image of the Baxter measure projected on dense $O(n)$ loop configurations has a density proportional to $\sqrt{q}^{-|\{{\rm non-contractible\ loops}\}|}$ wrt the random-cluster measure (in its loop representation).

\begin{figure}[htb]
\begin{center}
\scalebox{.5}{\includegraphics{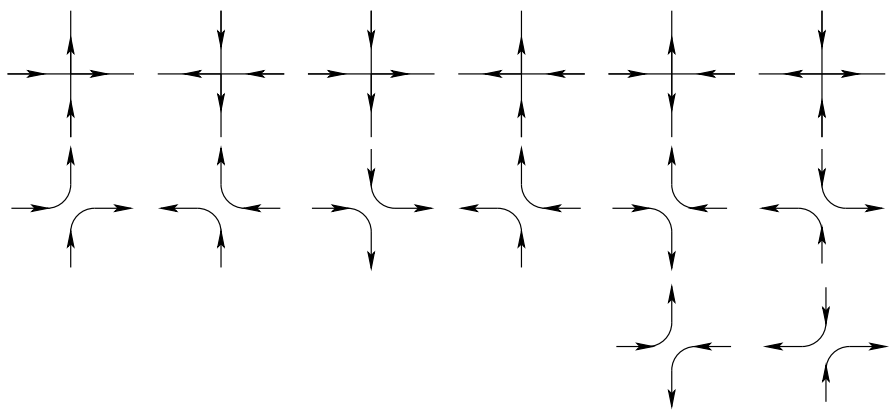}}
\end{center}
\caption{6V configurations and corresponding plaquette configurations for oriented loops (on each column)}
\label{Fig:6VBaxter}
\end{figure}

On the other hand, if one retains orientations (at each edge midpoint of $\dmd$) and forget loops, one obtains a six-vertex configuration on $\dmd^\dg$ (see Figure \ref{Fig:6VBaxter}), with weights:
$$\omega_1=\cdots=\omega_4=1, \omega_5=\omega_6=e^{i\pi s}+e^{-i\pi s}$$
so that $a=b=1$, $c^2=2+\sqrt q$.

It is worth pointing out that in this mapping, the parafermionic FK observables are mapped to magnetic-electric six-vertex observables. 
More generally, it is interesting to try to understand which natural observables may be expressed in both representations. Let us discuss topological observables (\cite{DiFSalZubcoulomb,Pinsontorus,ArgPotts}).

Fix $\omega=\alpha dx+\beta dy$ (thinking of the graph $\Gamma$  as having small mesh $\delta$ and being embedded on a macroscopic torus $\Sigma=\C/(\Z+\tau\Z)$, $\alpha,\beta$ fixed, and $\tau\in i(0,\infty)$ for simplicity; $A,B$ are the standard cycles on $\Sigma$). Let us modify the local Baxter weights as follows: the weight of an oriented quarter circle QT is $\exp(is\frac\pi 2)\exp(i\int_{QT}\omega)$. The weight of an oriented loop is now $\exp(\pm 2i\pi s)$ for a contractible loop and $\exp(i\int_L\omega)$ for a non-contractible loop. The possibilities for non-contractible loops are: no such loops (in which case either a cluster or a dual cluster contains paths homotopic to $A$ and $B$); 
$2k$ loops homotopic to $mA+nB$, $\gcd(m,n)=1$ (in which case there are $k$ clusters homotopic to a simple path representing $mA+nB$ separated by $k$ dual clusters). Let us denote $p_{k,m,n}$ the probability of this event under the random cluster measure. In this case, summing over orientations yields a relative weight $(2\cos(m\alpha+n\beta))^{2k}$.

On the six-vertex side, each north (resp. west, south, east) pointing edge gets a relative weight $e^{i\beta\delta}$ (resp. $e^{-i\alpha\delta},e^{-i\beta\delta},e^{i\alpha\delta}$). Collecting the weights of edges along an horizontal (resp. vertical) path on $\Gamma^\dg$ yields $\exp(-i\frac{2\beta\delta}\pi\int_A J)$ (resp. $\exp(i\frac{2\alpha\delta}\pi\int_B J)$).

Summarizing, setting 
$$f(\alpha,\beta)=p_0+\frac 12\sum_{k\geq 1}\sum_{\gcd(m,n)=1}p_{k,m,n}(2\cos(m\alpha+n\beta)q^{-1/2})^{2k}$$
we have the exact identity:
$$\frac{f(\alpha,\beta)}{f(0,0)}=\E_{6V}\left(\exp\left(2i\frac{\alpha}\pi\int_B J-2i\frac{\beta\Im\tau}\pi\int_A J\right)\right)$$
Admitting that the $6V$ height function converges to a compactified free field expresses the RHS as a theta function.\\

{\bf The free fermion point.}\\

\begin{figure}[htb]
\begin{center}
\scalebox{.5}{\includegraphics{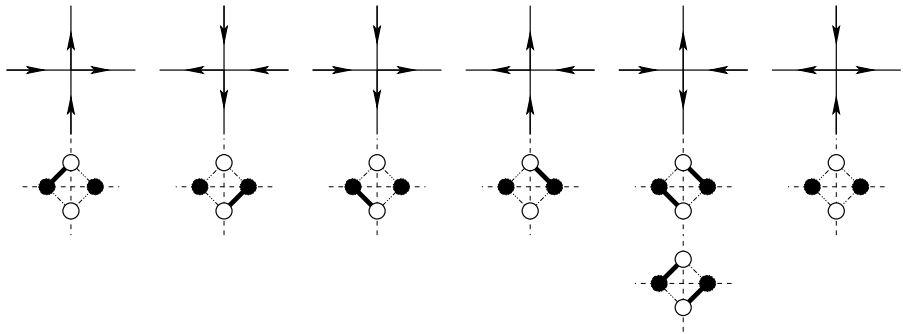}}
\end{center}
\caption{6V configurations and corresponding local dimer configuration (bold edge: dimer)}
\label{Fig:6Vdimer}
\end{figure}

For $c^2=a^2+b^2$, the (square) six-vertex model may be mapped exactly on a (weighted) dimer model on its medial graph. Fan and Wu (\cite{FanWu}) first showed the existence of a such a dimer representation in the more general context of the eight-vertex model. The mapping presented here is in particular apparent in the study of the six-vertex model with ``domain wall boundary conditions" (\cite{ICKDWBC}) and its relation with alternating sign matrices and the ``aztec diamond" problem for dimers (\cite{EKLP1,EKLP2,FerSpo6v}). The correspondence between dimer and 6V configurations is represented in Figure \ref{Fig:6Vdimer}.
Specifically, we start from $\Gamma$ which is modelled on the square lattice and carries a 6V configuration. Its (planar) derived graph $D\Gamma$ is also modelled on the square lattice; it has a 2-coloring, with black vertices corresponding to horizontal edges of $\Gamma$ and white vertices corresponding to vertical edges of $\Gamma$.
A dimer configuration (or {\em perfect matching}) on $D\Gamma$ is a subset ${\mf m}$ of edges of $D\Gamma$ such that each vertex of $D\Gamma$ is an endpoint of exactly one edge in ${\mf m}$.

In the 6V$\leftrightarrow$ dimer correspondence, if $(bw)\in{\mf m}$, with $b$ (resp. $w$) the midpoint of an horizontal (resp. vertical) edge of $\Gamma$, then: $\overrightarrow{bw}$ points in the direction of the edge through $b$ (with its 6V orientation) $\pm\frac\pi 4$; $\overrightarrow{wb}$ points in the opposite direction of the edge through $w$ (with its 6V orientation) $\pm\frac\pi 4$. It is easy to see that this prescription maps dimer tilings to 6V configuration, and that the tiling may be reconstructed from the 6V configuration up to local ambiguities (see type 5 in Figure \ref{Fig:6Vdimer}). One may also introduce defects: for instance, inserting edge defects at $b,w$ (with respective magnetic charges $-\frac 12,\frac 12$) corresponds exactly to considering dimer configurations with {\em monomer defects} at $b,w$ (\cite{FisSte2}).

In the dimer model (\cite{KenIAS}), each dimer configuration receives a weight
$$w({\mf m})=\prod_{e\in{\mf m}}w_e$$
where the $w_e$'s are fixed edge weights. Here, set $w_e=\cos(\theta)$ for $SW-NE$ edges and $w_e=\sin(\theta)$ for $SE-NW$ edges. Plainly, the corresponding 6V configuration has vertex weights:
$$(a:b:c)=(\cos(\theta):\sin(\theta):1)$$
i.e. we recover the {\em free fermion point} of the six-vertex model: $a^2+b^2=c^2$ or $\Delta=0$. 

In this case, we see that, up to normalization, the dimer height function (as defined in \cite{ThuCon,KenIAS} for the dimer model on a bipartite planar graph, taken here on every other face) is the height function of the corresponding 6V configurations. In this case, building on Kasteleyn's determinantal solution of the dimer model (\cite{Kassquare}), one can verify that the thermodynamic limit is described in some details by a (compactified) free field with coupling constant $g=2$. In particular, one can handle infinite volume Gibbs measures (\cite{CKP,KOS}); current correlations and scalar field limit (\cite{KendominoGFF}; for current correlations, it is convenient to consider height variations on the even sublattice of $\Gamma^\dg$); compactified free field limit (for an underlying periodic graph, \cite{Dubtors}); magnetic and electric correlators (\cite{Dubtors}).\\

{\bf Acknowledgments.}  These notes are based on lectures given at the 2011 Cornell Probability Summer School; I wish to thank its organizers for their kind hospitality.


\begin{thebibliography}{99}

\bibitem{ArgPotts}
\begin{barticle}[mr]
\bauthor{\bsnm{Arguin},~\bfnm{Louis-Pierre}\binits{L.-P.}}
(\byear{2002}).
\btitle{Homology of {F}ortuin-{K}asteleyn clusters of {P}otts models on the
  torus}.
\bjournal{J. Stat. Phys.}
\bvolume{109}
\bpages{301--310}.
\bid{doi={10.1023/A:1019979326380}, issn={0022-4715}, mr={1927924}}
\bptok{imsref}%
\end{barticle}
\endbibitem

\bibitem{Baxterexact}
\begin{bbook}[mr]
\bauthor{\bsnm{Baxter},~\bfnm{Rodney~J.}\binits{R.~J.}}
(\byear{1989}).
\btitle{Exactly Solved Models in Statistical Mechanics}.
\bpublisher{Academic Press [Harcourt Brace Jovanovich Publishers]},
  \baddress{London}.
\bnote{Reprint of the 1982 original}.
\bid{mr={0998375}}
\bptok{imsref}%
\end{bbook}
\endbibitem

\bibitem{CKP}
\begin{barticle}[mr]
\bauthor{\bsnm{Cohn},~\bfnm{Henry}\binits{H.}},
  \bauthor{\bsnm{Kenyon},~\bfnm{Richard}\binits{R.}} \AND
  \bauthor{\bsnm{Propp},~\bfnm{James}\binits{J.}}
(\byear{2001}).
\btitle{A variational principle for domino tilings}.
\bjournal{J. Amer. Math. Soc.}
\bvolume{14}
\bpages{297--346 (electronic)}.
\bid{doi={10.1090/S0894-0347-00-00355-6}, issn={0894-0347}, mr={1815214}}
\bptok{imsref}%
\end{barticle}
\endbibitem

\bibitem{DiF}
\begin{bbook}[mr]
\bauthor{\bsnm{Di~Francesco},~\bfnm{Philippe}\binits{P.}},
  \bauthor{\bsnm{Mathieu},~\bfnm{Pierre}\binits{P.}} \AND
  \bauthor{\bsnm{S{\'e}n{\'e}chal},~\bfnm{David}\binits{D.}}
(\byear{1997}).
\btitle{Conformal Field Theory}.
\bseries{Graduate Texts in Contemporary Physics}.
\bpublisher{Springer}, \baddress{New York}.
\bid{mr={1424041}}
\bptok{imsref}%
\end{bbook}
\endbibitem

\bibitem{DiFSalZubcoulomb}
\begin{barticle}[mr]
\bauthor{\bsnm{Di~Francesco},~\bfnm{P.}\binits{P.}},
  \bauthor{\bsnm{Saleur},~\bfnm{H.}\binits{H.}} \AND
  \bauthor{\bsnm{Zuber},~\bfnm{J.~B.}\binits{J.~B.}}
(\byear{1987}).
\btitle{Relations between the {C}oulomb gas picture and conformal invariance of
  two-dimensional critical models}.
\bjournal{J. Stat. Phys.}
\bvolume{49}
\bpages{57--79}.
\bid{issn={0022-4715}, mr={0923852}}
\bptok{imsref}%
\end{barticle}
\endbibitem

\bibitem{Dubtors}
\begin{bmisc}[auto:STB|2011/10/17|13:52:43]
\bauthor{\bsnm{{Dub{\'e}dat}},~\bfnm{J.}\binits{J.}}
\bhowpublished{{Dimers and analytic torsion I}. {\it arXiv:1110.2808}, 2011}.
\bptok{imsref}%
\end{bmisc}
\endbibitem

\bibitem{EdwSok}
\begin{barticle}[mr]
\bauthor{\bsnm{Edwards},~\bfnm{Robert~G.}\binits{R.~G.}} \AND
  \bauthor{\bsnm{Sokal},~\bfnm{Alan~D.}\binits{A.~D.}}
(\byear{1988}).
\btitle{Generalization of the {F}ortuin-{K}asteleyn-{S}wendsen-{W}ang
  representation and {M}onte {C}arlo algorithm}.
\bjournal{Phys. Rev. D (3)}
\bvolume{38}
\bpages{2009--2012}.
\bid{doi={10.1103/PhysRevD.38.2009}, issn={0556-2821}, mr={0965465}}
\bptok{imsref}%
\end{barticle}
\endbibitem

\bibitem{EKLP1}
\begin{barticle}[mr]
\bauthor{\bsnm{Elkies},~\bfnm{Noam}\binits{N.}},
  \bauthor{\bsnm{Kuperberg},~\bfnm{Greg}\binits{G.}},
  \bauthor{\bsnm{Larsen},~\bfnm{Michael}\binits{M.}} \AND
  \bauthor{\bsnm{Propp},~\bfnm{James}\binits{J.}}
(\byear{1992}).
\btitle{Alternating-sign matrices and domino tilings. {I}}.
\bjournal{J. Algebraic Combin.}
\bvolume{1}
\bpages{111--132}.
\bid{doi={10.1023/A:1022420103267}, issn={0925-9899}, mr={1226347}}
\bptok{imsref}%
\end{barticle}
\endbibitem

\bibitem{EKLP2}
\begin{barticle}[mr]
\bauthor{\bsnm{Elkies},~\bfnm{Noam}\binits{N.}},
  \bauthor{\bsnm{Kuperberg},~\bfnm{Greg}\binits{G.}},
  \bauthor{\bsnm{Larsen},~\bfnm{Michael}\binits{M.}} \AND
  \bauthor{\bsnm{Propp},~\bfnm{James}\binits{J.}}
(\byear{1992}).
\btitle{Alternating-sign matrices and domino tilings. {II}}.
\bjournal{J. Algebraic Combin.}
\bvolume{1}
\bpages{219--234}.
\bid{doi={10.1023/A:1022483817303}, issn={0925-9899}, mr={1194076}}
\bptok{imsref}%
\end{barticle}
\endbibitem

\bibitem{FanWu}
\begin{barticle}[auto:STB|2011/10/17|13:52:43]
\bauthor{\bsnm{Fan},~\bfnm{C.}\binits{C.}} \AND
  \bauthor{\bsnm{Wu},~\bfnm{F.~Y.}\binits{F.~Y.}}
(\byear{Aug 1970}).
\btitle{General lattice model of phase transitions}.
\bjournal{Phys. Rev. B}
\bvolume{2}
\bpages{723--733}.
\bptok{imsref}%
\end{barticle}
\endbibitem

\bibitem{FZZN}
\begin{barticle}[mr]
\bauthor{\bsnm{Fateev},~\bfnm{V.~A.}\binits{V.~A.}} \AND
  \bauthor{\bsnm{Zamolodchikov},~\bfnm{A.~B.}\binits{A.~B.}}
(\byear{1982}).
\btitle{Self-dual solutions of the star-triangle relations in
  {$Z\sb{N}$}-models}.
\bjournal{Phys. Lett. A}
\bvolume{92}
\bpages{37--39}.
\bid{doi={10.1016/0375-9601(82)90736-8}, issn={0375-9601}, mr={0677808}}
\bptok{imsref}%
\end{barticle}
\endbibitem

\bibitem{FerSpo6v}
\begin{barticle}[mr]
\bauthor{\bsnm{Ferrari},~\bfnm{Patrik~L.}\binits{P.~L.}} \AND
  \bauthor{\bsnm{Spohn},~\bfnm{Herbert}\binits{H.}}
(\byear{2006}).
\btitle{Domino tilings and the six-vertex model at its free-fermion point}.
\bjournal{J. Phys. A}
\bvolume{39}
\bpages{10297--10306}.
\bid{doi={10.1088/0305-4470/39/33/003}, issn={0305-4470}, mr={2256593}}
\bptok{imsref}%
\end{barticle}
\endbibitem

\bibitem{FisSte2}
\begin{barticle}[mr]
\bauthor{\bsnm{Fisher},~\bfnm{Michael~E.}\binits{M.~E.}} \AND
  \bauthor{\bsnm{Stephenson},~\bfnm{John}\binits{J.}}
(\byear{1963}).
\btitle{Statistical mechanics of dimers on a plane lattice. {II}. {D}imer
  correlations and monomers}.
\bjournal{Phys. Rev. (2)}
\bvolume{132}
\bpages{1411--1431}.
\bid{mr={0158705}}
\bptok{imsref}%
\end{barticle}
\endbibitem

\bibitem{GawCFT}
\begin{bincollection}[mr]
\bauthor{\bsnm{Gaw{\c{e}}dzki},~\bfnm{Krzysztof}\binits{K.}}
(\byear{1999}).
\btitle{Lectures on conformal field theory}.
In \bbooktitle{Quantum Fields and Strings: A Course for Mathematicians, {V}ol.
  1, 2 ({P}rinceton, {NJ}, 1996/1997)}
\bpages{727--805}.
\bpublisher{Amer. Math. Soc.}, \baddress{Providence, RI}.
\bid{mr={1701610}}
\bptok{imsref}%
\end{bincollection}
\endbibitem

\bibitem{GeoGibbs}
\begin{bbook}[mr]
\bauthor{\bsnm{Georgii},~\bfnm{Hans-Otto}\binits{H.-O.}}
(\byear{2011}).
\btitle{Gibbs Measures and Phase Transitions},
\bedition{2nd} ed.
\bseries{de Gruyter Studies in Mathematics}
\bvolume{9}.
\bpublisher{de Gruyter}, \baddress{Berlin}.
\bid{mr={2807681}}
\bptok{imsref}%
\end{bbook}
\endbibitem

\bibitem{GlimmJaffe}
\begin{bbook}[mr]
\bauthor{\bsnm{Glimm},~\bfnm{James}\binits{J.}} \AND
  \bauthor{\bsnm{Jaffe},~\bfnm{Arthur}\binits{A.}}
(\byear{1987}).
\btitle{Quantum Physics}, \bedition{2nd} ed.
\bpublisher{Springer}, \baddress{New York}.
\bnote{A functional integral point of view}.
\bid{mr={0887102}}
\bptok{imsref}%
\end{bbook}
\endbibitem

\bibitem{Grimperc}
\begin{bbook}[mr]
\bauthor{\bsnm{Grimmett},~\bfnm{Geoffrey}\binits{G.}}
(\byear{1999}).
\btitle{Percolation},
\bedition{2nd} ed.
\bseries{Grundlehren der Mathematischen Wissenschaften [Fundamental Principles
  of Mathematical Sciences]}
\bvolume{321}.
\bpublisher{Springer}, \baddress{Berlin}.
\bid{mr={1707339}}
\bptok{imsref}%
\end{bbook}
\endbibitem

\bibitem{GrimFK}
\begin{bbook}[mr]
\bauthor{\bsnm{Grimmett},~\bfnm{Geoffrey}\binits{G.}}
(\byear{2006}).
\btitle{The Random-cluster Model}.
\bseries{Grundlehren der Mathematischen Wissenschaften [Fundamental Principles
  of Mathematical Sciences]}
\bvolume{333}.
\bpublisher{Springer}, \baddress{Berlin}.
\bid{doi={10.1007/978-3-540-32891-9}, mr={2243761}}
\bptok{imsref}%
\end{bbook}
\endbibitem

\bibitem{Grimgraphs}
\begin{bbook}[mr]
\bauthor{\bsnm{Grimmett},~\bfnm{Geoffrey}\binits{G.}}
(\byear{2010}).
\btitle{Probability on Graphs}.
\bseries{Institute of Mathematical Statistics Textbooks}
\bvolume{1}.
\bpublisher{Cambridge Univ. Press}, \baddress{Cambridge}.
\bnote{Random processes on graphs and lattices}.
\bid{mr={2723356}}
\bptok{imsref}%
\end{bbook}
\endbibitem

\bibitem{GrunbaumDFT}
\begin{barticle}[mr]
\bauthor{\bsnm{Gr{\"u}nbaum},~\bfnm{F.~Alberto}\binits{F.~A.}}
(\byear{1982}).
\btitle{The eigenvectors of the discrete {F}ourier transform: A version of the
  {H}ermite functions}.
\bjournal{J. Math. Anal. Appl.}
\bvolume{88}
\bpages{355--363}.
\bid{doi={10.1016/0022-247X(82)90199-8}, issn={0022-247X}, mr={0667064}}
\bptok{imsref}%
\end{barticle}
\endbibitem

\bibitem{IkhCar}
\begin{barticle}[mr]
\bauthor{\bsnm{Ikhlef},~\bfnm{Yacine}\binits{Y.}} \AND
  \bauthor{\bsnm{Cardy},~\bfnm{John}\binits{J.}}
(\byear{2009}).
\btitle{Discretely holomorphic parafermions and integrable loop models}.
\bjournal{J. Phys. A}
\bvolume{42}
\bpages{102001, 11}.
\bid{doi={10.1088/1751-8113/42/10/102001}, issn={1751-8113}, mr={2485852}}
\bptok{imsref}%
\end{barticle}
\endbibitem

\bibitem{ICKDWBC}
\begin{barticle}[mr]
\bauthor{\bsnm{Izergin},~\bfnm{A.~G.}\binits{A.~G.}},
  \bauthor{\bsnm{Coker},~\bfnm{D.~A.}\binits{D.~A.}} \AND
  \bauthor{\bsnm{Korepin},~\bfnm{V.~E.}\binits{V.~E.}}
(\byear{1992}).
\btitle{Determinant formula for the six-vertex model}.
\bjournal{J. Phys. A}
\bvolume{25}
\bpages{4315--4334}.
\bid{issn={0305-4470}, mr={1181591}}
\bptok{imsref}%
\end{barticle}
\endbibitem

\bibitem{KCdisorder}
\begin{barticle}[mr]
\bauthor{\bsnm{Kadanoff},~\bfnm{Leo~P.}\binits{L.~P.}} \AND
  \bauthor{\bsnm{Ceva},~\bfnm{Horacio}\binits{H.}}
(\byear{1971}).
\btitle{Determination of an operator algebra for the two-dimensional {I}sing
  model}.
\bjournal{Phys. Rev. B (3)}
\bvolume{3}
\bpages{3918--3939}.
\bid{mr={0389111}}
\bptok{imsref}%
\end{barticle}
\endbibitem

\bibitem{Kassquare}
\begin{barticle}[auto:STB|2011/10/17|13:52:43]
\bauthor{\bsnm{Kasteleyn},~\bfnm{P.~W.}\binits{P.~W.}}
(\byear{1961}).
\btitle{The statistics of dimers on a lattice. i. the number of dimer
  arrangements on a quadratic lattice}.
\bjournal{Physica}
\bvolume{27}
\bpages{1209--1225}.
\bptok{imsref}%
\end{barticle}
\endbibitem

\bibitem{Kendouble}
\begin{bmisc}[auto:STB|2011/10/17|13:52:43]
\bauthor{\bsnm{Kenyon},~\bfnm{R.}\binits{R.}}
\bhowpublished{Conformal invariance of loops in the double-dimer model. {\it
  preprint, arXiv:1105.4158}, 2011}.
\bptok{imsref}%
\end{bmisc}
\endbibitem

\bibitem{KendominoGFF}
\begin{barticle}[mr]
\bauthor{\bsnm{Kenyon},~\bfnm{Richard}\binits{R.}}
(\byear{2001}).
\btitle{Dominos and the {G}aussian free field}.
\bjournal{Ann. Probab.}
\bvolume{29}
\bpages{1128--1137}.
\bid{doi={10.1214/aop/1015345599}, issn={0091-1798}, mr={1872739}}
\bptok{imsref}%
\end{barticle}
\endbibitem

\bibitem{KenIAS}
\begin{bincollection}[mr]
\bauthor{\bsnm{Kenyon},~\bfnm{Richard}\binits{R.}}
(\byear{2009}).
\btitle{Lectures on dimers}.
In \bbooktitle{Statistical Mechanics}.
\bseries{IAS/Park City Math. Ser.}
\bvolume{16}
\bpages{191--230}.
\bpublisher{Amer. Math. Soc.}, \baddress{Providence, RI}.
\bid{mr={2523460}}
\bptok{imsref}%
\end{bincollection}
\endbibitem

\bibitem{KOS}
\begin{barticle}[mr]
\bauthor{\bsnm{Kenyon},~\bfnm{Richard}\binits{R.}},
  \bauthor{\bsnm{Okounkov},~\bfnm{Andrei}\binits{A.}} \AND
  \bauthor{\bsnm{Sheffield},~\bfnm{Scott}\binits{S.}}
(\byear{2006}).
\btitle{Dimers and amoebae}.
\bjournal{Ann. of Math. (2)}
\bvolume{163}
\bpages{1019--1056}.
\bid{doi={10.4007/annals.2006.163.1019}, issn={0003-486X}, mr={2215138}}
\bptok{imsref}%
\end{barticle}
\endbibitem

\bibitem{KWising1}
\begin{barticle}[mr]
\bauthor{\bsnm{Kramers},~\bfnm{H.~A.}\binits{H.~A.}} \AND
  \bauthor{\bsnm{Wannier},~\bfnm{G.~H.}\binits{G.~H.}}
(\byear{1941}).
\btitle{Statistics of the two-dimensional ferromagnet. {I}}.
\bjournal{Phys. Rev. (2)}
\bvolume{60}
\bpages{252--262}.
\bid{mr={0004803}}
\bptok{imsref}%
\end{barticle}
\endbibitem

\bibitem{KWising2}
\begin{barticle}[mr]
\bauthor{\bsnm{Kramers},~\bfnm{H.~A.}\binits{H.~A.}} \AND
  \bauthor{\bsnm{Wannier},~\bfnm{G.~H.}\binits{G.~H.}}
(\byear{1941}).
\btitle{Statistics of the two-dimensional ferromagnet. {II}}.
\bjournal{Phys. Rev. (2)}
\bvolume{60}
\bpages{263--276}.
\bid{mr={0004804}}
\bptok{imsref}%
\end{barticle}
\endbibitem

\bibitem{Liebsquare}
\begin{barticle}[auto:STB|2011/10/17|13:52:43]
\bauthor{\bsnm{Lieb},~\bfnm{E.~H.}\binits{E.~H.}}
(\byear{Oct 1967}).
\btitle{Residual entropy of square ice}.
\bjournal{Phys. Rev.}
\bvolume{162}
\bpages{162--172}.
\bptok{imsref}%
\end{barticle}
\endbibitem

\bibitem{MWising}
\begin{bmisc}[auto:STB|2011/10/17|13:52:43]
\bauthor{\bsnm{McCoy},~\bfnm{B.}\binits{B.}} \AND
  \bauthor{\bsnm{Wu},~\bfnm{T.}\binits{T.}}
\bhowpublished{{\it The two-dimensional Ising model}. Harvard Univ. Press,
  Boston, MA, 1973}.
\bptok{imsref}%
\end{bmisc}
\endbibitem

\bibitem{Mercat}
\begin{barticle}[mr]
\bauthor{\bsnm{Mercat},~\bfnm{Christian}\binits{C.}}
(\byear{2001}).
\btitle{Discrete {R}iemann surfaces and the {I}sing model}.
\bjournal{Comm. Math. Phys.}
\bvolume{218}
\bpages{177--216}.
\bid{doi={10.1007/s002200000348}, issn={0010-3616}, mr={1824204}}
\bptok{imsref}%
\end{barticle}
\endbibitem

\bibitem{NienCoulomb}
\begin{barticle}[mr]
\bauthor{\bsnm{Nienhuis},~\bfnm{Bernard}\binits{B.}}
(\byear{1984}).
\btitle{Critical behavior of two-dimensional spin models and charge asymmetry
  in the {C}oulomb gas}.
\bjournal{J. Stat. Phys.}
\bvolume{34}
\bpages{731--761}.
\bid{issn={0022-4715}, mr={0751711}}
\bptok{imsref}%
\end{barticle}
\endbibitem

\bibitem{NieKnopotts}
\begin{barticle}[auto:STB|2011/10/17|13:52:43]
\bauthor{\bsnm{Nienhuis},~\bfnm{B.}\binits{B.}} \AND
  \bauthor{\bsnm{Knops},~\bfnm{H.~J.~F.}\binits{H.~J.~F.}}
(\byear{1872--1875, Aug}).
\btitle{Spinor exponents for the two-dimensional potts model}.
\bjournal{Phys. Rev. B}
\bvolume{32}
\bpages{1985}.
\bptok{imsref}%
\end{barticle}
\endbibitem

\bibitem{Palmerplanar}
\begin{bbook}[mr]
\bauthor{\bsnm{Palmer},~\bfnm{John}\binits{J.}}
(\byear{2007}).
\btitle{Planar {I}sing Correlations}.
\bseries{Progress in Mathematical Physics}
\bvolume{49}.
\bpublisher{Birkh\"auser}, \baddress{Boston, MA}.
\bid{mr={2332010}}
\bptok{imsref}%
\end{bbook}
\endbibitem

\bibitem{Pinsontorus}
\begin{barticle}[mr]
\bauthor{\bsnm{Pinson},~\bfnm{Haru~T.}\binits{H.~T.}}
(\byear{1994}).
\btitle{Critical percolation on the torus}.
\bjournal{J. Stat. Phys.}
\bvolume{75}
\bpages{1167--1177}.
\bid{issn={0022-4715}, mr={1285297}}
\bptok{imsref}%
\end{barticle}
\endbibitem

\bibitem{Resh6V}
\begin{bincollection}[mr]
\bauthor{\bsnm{Reshetikhin},~\bfnm{N.}\binits{N.}}
(\byear{2010}).
\btitle{Lectures on the integrability of the six-vertex model}.
In \bbooktitle{Exact Methods in Low-dimensional Statistical Physics and Quantum
  Computing}
\bpages{197--266}.
\bpublisher{Oxford Univ. Press}, \baddress{Oxford}.
\bid{mr={2668647}}
\bptok{imsref}%
\end{bincollection}
\endbibitem

\bibitem{CarRiv}
\begin{barticle}[mr]
\bauthor{\bsnm{Riva},~\bfnm{V.}\binits{V.}} \AND
  \bauthor{\bsnm{Cardy},~\bfnm{J.}\binits{J.}}
(\byear{2006}).
\btitle{Holomorphic parafermions in the {P}otts model and stochastic {L}oewner
  evolution}.
\bjournal{J. Stat. Mech. Theory Exp.}
\bvolume{12}
\bpages{P12001, 19 pp. (electronic)}.
\bid{issn={1742-5468}, mr={2280251}}
\bptok{imsref}%
\end{barticle}
\endbibitem

\bibitem{Rudinfourier}
\begin{bbook}[mr]
\bauthor{\bsnm{Rudin},~\bfnm{Walter}\binits{W.}}
(\byear{1990}).
\btitle{Fourier Analysis on Groups}.
\bseries{Wiley Classics Library}.
\bpublisher{Wiley}, \baddress{New York}.
\bnote{Reprint of the 1962 original, A Wiley-Interscience Publication}.
\bid{mr={1038803}}
\bptok{imsref}%
\end{bbook}
\endbibitem

\bibitem{Savitdual}
\begin{barticle}[mr]
\bauthor{\bsnm{Savit},~\bfnm{Robert}\binits{R.}}
(\byear{1982}).
\btitle{Duality transformations for general abelian systems}.
\bjournal{Nuclear Phys. B}
\bvolume{200}
\bpages{233--248}.
\bid{doi={10.1016/0550-3213(82)90085-2}, issn={0550-3213}, mr={0643588}}
\bptok{imsref}%
\end{barticle}
\endbibitem

\bibitem{SimonPphi}
\begin{bbook}[mr]
\bauthor{\bsnm{Simon},~\bfnm{Barry}\binits{B.}}
(\byear{1974}).
\btitle{The {$P(\phi )\sb{2}$} {E}uclidean (quantum) Field Theory}.
\bpublisher{Princeton Univ. Press}, \baddress{Princeton, N.J.}
\bnote{Princeton Series in Physics}.
\bid{mr={0489552}}
\bptok{imsref}%
\end{bbook}
\endbibitem

\bibitem{SmiICM}
\begin{bincollection}[mr]
\bauthor{\bsnm{Smirnov},~\bfnm{Stanislav}\binits{S.}}
(\byear{2006}).
\btitle{Towards conformal invariance of 2{D} lattice models}.
In \bbooktitle{International {C}ongress of {M}athematicians. {V}ol. {II}}
\bpages{1421--1451}.
\bpublisher{Eur. Math. Soc., Z\"urich}.
\bid{mr={2275653}}
\bptok{imsref}%
\end{bincollection}
\endbibitem

\bibitem{Smiising}
\begin{barticle}[mr]
\bauthor{\bsnm{Smirnov},~\bfnm{Stanislav}\binits{S.}}
(\byear{2010}).
\btitle{Conformal invariance in random cluster models. {I}. {H}olomorphic
  fermions in the {I}sing model}.
\bjournal{Ann. of Math. (2)}
\bvolume{172}
\bpages{1435--1467}.
\bid{doi={10.4007/annals.2010.172.1441}, issn={0003-486X}, mr={2680496}}
\bptok{imsref}%
\end{barticle}
\endbibitem

\bibitem{ThuCon}
\begin{barticle}[mr]
\bauthor{\bsnm{Thurston},~\bfnm{William~P.}\binits{W.~P.}}
(\byear{1990}).
\btitle{Conway's tiling groups}.
\bjournal{Amer. Math. Monthly}
\bvolume{97}
\bpages{757--773}.
\bid{doi={10.2307/2324578}, issn={0002-9890}, mr={1072815}}
\bptok{imsref}%
\end{barticle}
\endbibitem

\bibitem{VB6VSOS}
\begin{barticle}[auto:STB|2011/10/17|13:52:43]
\bauthor{\bparticle{van} \bsnm{Beijeren},~\bfnm{H.}\binits{H.}}
(\byear{May 1977}).
\btitle{Exactly solvable model for the roughening transition of a crystal
  surface}.
\bjournal{Phys. Rev. Lett.}
\bvolume{38}
\bpages{993--996}.
\bptok{imsref}%
\end{barticle}
\endbibitem

\bibitem{W1}
\begin{bincollection}[mr]
\bauthor{\bsnm{Werner},~\bfnm{Wendelin}\binits{W.}}
(\byear{2004}).
\btitle{Random planar curves and {S}chramm-{L}oewner evolutions}.
In \bbooktitle{Lectures on Probability Theory and Statistics}.
\bseries{Lecture Notes in Math.}
\bvolume{1840}
\bpages{107--195}.
\bpublisher{Springer}, \baddress{Berlin}.
\bid{mr={2079672}}
\bptok{imsref}%
\end{bincollection}
\endbibitem

\bibitem{WWduality}
\begin{barticle}[auto:STB|2011/10/17|13:52:43]
\bauthor{\bsnm{Wu},~\bfnm{F.~Y.}\binits{F.~Y.}} \AND
  \bauthor{\bsnm{Wang},~\bfnm{Y.~K.}\binits{Y.~K.}}
(\byear{1976}).
\btitle{Duality transformation in a many-component spin model}.
\bjournal{J. Math. Phys.}
\bvolume{17}
\bpages{439--440}.
\bptok{imsref}%
\end{barticle}
\endbibitem

\end{thebibliography}
\end{document}